\documentclass{llncs}

\usepackage[utf8]{inputenc}
\usepackage[T1]{fontenc}
\usepackage{enumerate}
\usepackage{graphicx}
\usepackage{amsfonts}
\usepackage{amssymb}
\usepackage{amsmath}
\usepackage{wasysym}
\usepackage{hyperref}
\usepackage{mathtools}
\usepackage{bibentry}
\usepackage{comment}
\usepackage{paralist}
\usepackage{enumitem}
\usepackage{todonotes}
\usepackage{hyperref}
\usepackage{thmtools}
\usepackage{float}

\newenvironment{observation}{{\bf Observation }}{}

\graphicspath{{./figures/}}%helpful if your graphic files are in another directory

\def\computationproblem#1#2#3{
  \begin{center}
  \begin{tabular}{rp{0.8\textwidth}}
  {\sc Problem:\enspace}&#1\\
  {\sc Input:\enspace}&#2\\
  {\sc Question:\enspace}&#3\\
  \end{tabular}
  \end{center}
}

\title{Graph covers and semi-covers: Who is stronger?\thanks{An extended abstract of this paper has been presented at Eurocomb'23~\cite{eurocomb23}.}}

\author{Jan Kratochvíl\inst{1}\orcidID{0000-0002-2620-6133} 
\and Roman Nedela\inst{2}\orcidID{0000-0002-9826-704X}
}

\institute{
Department of Applied Mathematics, Faculty of Mathematics and Physics, Charles University, Prague, Czech Republic, \url{honza@kam.mff.cuni.cz}
\and
Faculty of Applied Sciences, University of West Bohemia, Pilsen, Czech Republic, \url{nedela@savbb.sk}
}
\begin{document}

\pagestyle{plain}

\maketitle

\begin{abstract}
The notion of graph cover, also known as locally bijective homomorphism, is a discretization of covering spaces known from general topology. It is a pair of incidence-preserving vertex- and edge-mappings between two graphs, the edge-component being bijective on the edge-neighborhoods of every vertex and its image. In line with the current trends in topological graph theory and its applications in mathematical physics, graphs are considered in the most relaxed form and as such they may contain multiple edges, loops and semi-edges.

Nevertheless, simple graphs (binary structures without multiple edges, loops, or semi-edges) play an important role. It has been conjectured in [Bok et al.: List covering of regular multigraphs, Proceedings IWOCA 2022, LNCS 13270, pp. 228--242] that for every fixed graph $H$, deciding if a graph covers $H$ is either polynomial time solvable for arbitrary input graphs, or NP-complete for simple ones. A graph $A$ is called stronger than a graph $B$ if every simple graph that covers $A$ also covers $B$. This notion was defined and found useful for NP-hardness reductions for disconnected graphs in [Bok et al.: Computational complexity of covering disconnected multigraphs, Proceedings FCT 2022, LNCS 12867, pp. 85--99]. It was conjectured in [Kratochv\'{\i}l: Towards strong dichotomy of graphs covers, GROW 2022 - Book of open problems, p. 10, {\tt https://grow.famnit.upr.si/GROW-BOP.pdf}] that if $A$ has no semi-edges, then $A$ is stronger than $B$ if and only if $A$ covers $B$. We prove this conjecture for cubic one-vertex graphs, and we also justify it for all cubic graphs $A$ with at most 4 vertices.        
\end{abstract}

\vskip1cm

%\todo{ORCID Nedela zkontrolovat - RN}
%\todo{abstract - JK}
%\todo{uvod - obecne nakryti - RN}
%\todo{uvod - motivace pro trojuhelnicek a slozitost - JK}
%\todo{uvod - nase vysledky - JK bude pokracovat RN}
%\todo{presunout hranove barveni odnekud nekam - JK}

%$\leadsto $

%$\longrightarrow$

\section{Introduction}\label{sec:Intro}

Combinatorial treatment of graph coverings had its
primary incentive in the solution of Heawood's Map Colour Problem
due to Ringel, Youngs and others \cite{Ringel1968,Ringel2012}. That coverings
underlie the techniques that led to the eventual solution of the
problem was recognized by Alpert and Gross \cite{gross1974}. These ideas
further crystallized in 1974 in the work of Gross \cite{gross1974voltage} where
voltage graphs were introduced as a means of a purely combinatorial
description of regular graph coverings. In parallel, the very same
idea appeared in Biggs's monograph \cite{biggs1974}. Much of the theory
of combinatorial graph coverings in its own right was subsequently
developed by Gross and Tucker in the seventies. We refer the reader
to \cite{gross2001,white1985} and the references therein. 

Since the book of Gross and Tucker \cite{gross2001}
was published, it remains the most popular textbook as well as
the reference book covering the topic.
In \cite{malnic2000}  the combinatorial theory of graph coverings
and voltage assignments was established and extended onto a more general class of graphs which include edges with free ends (called semi-edges). The new concept of a graph proved to be useful in applications as well as in theoretical considerations.
From the theoretical point of view, graphs with semi-edges arise
as quotients of standard graphs by groups of automorphisms which are semiregular
on vertices and darts (arcs) but may fix edges. They may be viewed
as 1-dimensional analogues of 2- and 3-dimensional orbifolds.
The fruitful concept of an orbifold comes from studies
of geometry of $n$-manifolds. Orbifolds and related concepts are
implicitly included in the work of pioneers such as Henri Poincare, or
Paul K\" obe. The first formal definition of an orbifold-like object was given by Ichiro Satake in 1956 in the context of Riemannian geometry. William Thurston, later in the mid 1970s defined and named the more general notion of orbifold as part of his study of hyperbolic structures. It is not a surprise that the concept of graphs with semi-edges was discovered in context of theoretical physics, as well 
(see \cite{getzler1998} for instance).

The study of the automorphism lifting problem in the context of
regular coverings of graphs had its main motivation in constructing
infinite families of highly transitive graphs. The first notable
contribution along these lines appeared, incidentally, in 1974 in
Biggs' monograph \cite{biggs1974} and in a paper of Djokovi\'{c}
\cite{djokovic1974}. Whereas Biggs gave a combinatorial sufficient
condition for a lifted group to be a split extension, Djokovi\'{c}
found a criterion, in terms of the fundamental group, for a group of
automorphisms to lift at all. A decade later, several different
sources added further motivations for studying the lifting problem.
These include: counting isomorphism classes of coverings and, more
generally, graph bundles, as considered by Hofmeister (1991)
and Kwak and Lee (1990,1992); constructions of regular maps on
surfaces based on covering space techniques due to Archdeacon,
Gvozdjak, Nedela, Richter, \v Sir\'{a}\v n, \v Skoviera and Surowski; and construction of
transitive graphs with prescribed degree of symmetry, for instance
by Du, Malni\v c, Nedela, Maru\v si\v c, Scapellato, Seifter,
Trofimov and Waller. Lifting and/or
projecting techniques play a prominent role also in the study of
imprimitive graphs, cf. Gardiner and Praeger \cite{gardiner1995} among
others. Nowadays the construction of a graph covering over a prescribed
graph is established as a useful technique allowing to construct effectively infinite families of graphs sharing prescribed properties. In particular, it was used to construct extremal regular graphs with fixed degree and diameter \cite{miller2012}, to construct cages and their approximations \cite{jajcay2011}, and in investigation of flows on graphs \cite{nedela2001} .

From the computational complexity point of view, Bodlaender~\cite{n:Bodlaender89} showed that deciding if a given graph $G$ covers a given graph $H$ (both graphs are part of the input) is NP-complete. Abello et al.~\cite{n:AFS91} asked for the complexity of this question when parameterized by the target graph $H$. They gave the first examples of graphs for which the problem, referred to as {\sc $H$-Cover}, is NP-complete or polynomial time solvable. It should be noted that in this seminal paper, both the parameter and the input graphs are allowed to have loops and multiple edges, but not semi-edges. The impact of semi-edges for the complexity issues is first discussed in Bok et al.~\cite{n:BFHJK-MFCS}. It is perhaps somewhat surprising that in all cases where the complexity of the {\sc $H$-Cover} problem is known to be NP-complete, it has been proved NP-complete even for simple graphs on input. This has been now conjectured to hold true in general, as the Strong Dichotomy Conjecture (cf. Conjecture~1 below) in~\cite{iwoca,BokFJKR24-algo}.

Bok et al.~\cite{BFJKS21FCT,BokFJKS24-dam} have proved that the Strong Dichotomy Conjecture holds true for all graphs, provided it holds true for connected ones. The curiosity of the NP-hardness reduction is its non-constructiveness. For two graphs $A$ and $B$, they use a simple graph $A'$ which covers $A$ and does not cover $B$, if such a simple graph exists (and an arbitrary simple cover of $A$ otherwise). However, there is no clue how to decide if such an $A'$ exists or not. (This paradox is not undermining the reduction, because the approach is used for fixed graphs $A$ and $B$, two of the connected components of the target graph $H$, to prove the existence of a polynomial time reduction between computational problems.) As a consequence, they defined a binary relation between connected graphs, saying that a graph $A$ is stronger than a graph $B$ if such a graph $A'$ does not exist. Our aim is to contribute to the study of this relation.

A snark is a simple 2-connected cubic graph which is not 3-edge-colorable. 
Interest in snarks was boosted by an observation by Heawood (1890) that the Four colour theorem is equivalent to the statement: there are no planar snarks. The most famous example of a snark is the Petersen graph. For about 80 years only few examples of non-trivial snarks were constructed until Isaacs (1975) introduced 
the infinite family of flower snarks and the operation of dot-product allowing
to construct a new nontrivial snark from two given ones. Recall that a snark is non-trivial if it is cyclically $4$-connected and of girth at least five.
Nowadays investigation of snarks is an active area of research due to the fact
that many long-standing conjectures on graphs (such as the 5-flow conjecture or the cycle double cover conjecture) can be reduced to problems on snarks, see \cite{watkins1989,fiol2018,MRS2022} and the references therein.

If we denote by $F(3,0)$ the one-vertex graph with 3 semi-edges, and by $F(1,1)$ the one-vertex graph with 1 semi-edge and 1 loop, a snark is a simple cubic graph that covers $F(1,1)$ but does not cover $F(3,0)$, i.e., a witness for the fact that $B=F(1,1)$ is not stronger than $A=F(3,0)$. By the Petersen theorem, a 2-connected cubic graph always contains a perfect matching          and hence it covers $F(1,1)$. In this sense every witness $A'$ for $A$ not being stronger than $B$ can be viewed as a generalized snark.

It is easy to see that  $A$ is stronger than $B$ whenever $A$ covers $B$. For all known pairs $A,B$ such that $A$ is stronger than $B$ and $A$ does not cover $B$, the graph $A$ contains semi-edges. In~\cite{n:GROW2022-open} it is conjectured that this is always the case (cf. Conjectures~2 and~3 below). In this paper, we justify these conjectures in several general situations, namely for cubic graphs.   

\section{Preliminaries}\label{sec:Prelim}

\subsection{Definitions}\label{subsec:Defs}

\begin{definition}
A {\em graph} is a finite set of vertices accompanied with a finite set of {\em edges}, where an edge is either a {\em loop}, or a {\em semi-edge}, or a {\em normal edge}. A normal edge is incident with 2 vertices and adds 1 to the degree of each of them. A loop is incident with a single vertex and adds 2 to its degree. A semi-edge is also incident with a single vertex, but adds only 1 to its degree. 
\end{definition}

As defined, we only consider undirected graphs.
However, graphs may have multiple loops and/or multiple semi-edges incident with the same vertex, and also multiple normal edges incident with the same pair of vertices. A graph is called {\em simple} if it has no loops, no semi-edges and no multiple normal edges.  The {\em edge-neighborhood} $E_G(u)$ of a vertex $u$ is the set of edges of $G$ incident with $u$.

\begin{figure}
\centering
{\includegraphics[width=\textwidth]{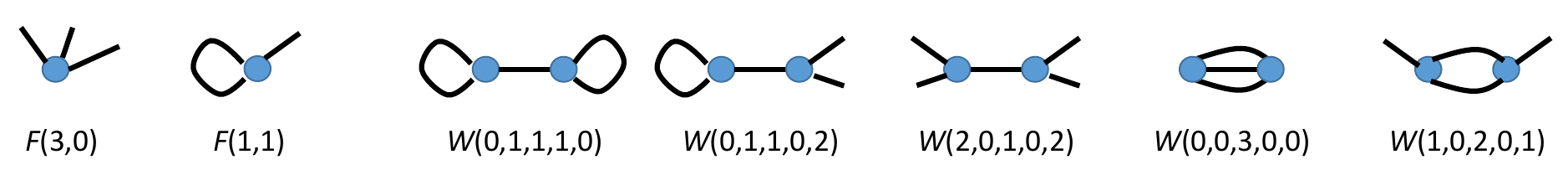}}
\caption{Connected cubic $F$ and $W$ graphs.}
\label{fig:examplesflowersanddumbells}
\end{figure}

\medskip\noindent
{\bf Notation.} The one-vertex graph with $a$ semi-edges and $b$ loops is denoted by $F(a,b)$. The two-vertex graph with $k$ semi-edges and $m$ loops incident with one if its vertices, $p$ loops and $q$ semi-edges incident with the other one, and $\ell$ normal edges incident with both of them is denoted by $W(k,m,\ell,p,q)$. See examples in Figure~\ref{fig:examplesflowersanddumbells}.

\begin{definition}
A {\em covering projection} from a graph $G$ to a connected graph $H$ is a pair of surjective mappings $f_V:V(G)\longrightarrow V(H)$ and $f_E:E(G)\longrightarrow E(H)$ such that

- $f_V$ is degree preserving,

- $f_E$ maps semi-edges onto semi-edges and loops onto loops, respectively, (normal edges may be mapped onto normal edges, loops, and semi-edges),

- $f_E$ is incidence preserving (i.e., if $e\in E(G)$ is incident with vertices $u,v\in V(G)$, then $f_E(e)$ is incident with $f_V(u)$ and $f_V(v)$, which may of course be the same vertex),

- $f_E$ is a local bijection on the edge-neighborhoods of any vertex and its image. %See~\cite{n:BFHJK-MFCS} for a more formal definition. 
In particular, the preimage of a normal edge incident with vertices $u,v\in V(H)$ is a disjoint union of normal edges, each incident with one vertex in $f_V^{-1}(u)$ and with one vertex in $f_V^{-1}(v)$, spanning $f_V^{-1}(u) \cup f_V^{-1}(v)$; the preimage of a loop incident with vertex $u\in V(H)$ is a disjoint union of cycles spanning $f_V^{-1}(u)$ (a loop is a cycle of length 1, and two parallel edges form a cycle of length 2); and the preimage of a semi-edge incident with a vertex $u\in V(H)$ is a disjoint union of semi-edges and normal edges spanning $f_V^{-1}(u)$. 
\end{definition}

If a graph $G$ allows a covering projection onto a graph $H$, we say that $G$ {\em covers} $H$, and we write $G\longrightarrow H$.

It is well known that in a covering projection to a connected graph, the sizes of preimages of all vertices of the target graph are the same. This implies that $|V(H)|$ divides $|V(G)|$ whenever $G\longrightarrow H$ for a connected graph $H$. We say that $G$ is a $k$-fold cover of $H$, with $k=\frac{|V(G)|}{|V(H)|}$ in such a case. 

We are interested in the following computational problem, parameterized by the target graph $H$.

\computationproblem{\sc $H$-Cover}{A graph $G$.}{Does $G$ cover $H$?}

Abello et al.~\cite{n:AFS91} raised the question of characterizing the complexity of {\sc $H$-Cover} for simple graphs $H$. Despite intensive effort and several general results, the complete characterization and even a conjecture on what are the easy and hard cases is not in sight. Bok et al.~\cite{n:BFHJK-MFCS} is the first paper that studies this question for (multi)graphs with semi-edges. A polynomial/NP-completeness dichotomy is believed in, and it has been conjectured in a stronger form in~\cite{iwoca}:

\medskip\noindent
{\bf Conjecture 1 (Strong Dichotomy Conjecture).} For every graph $H$, the {\sc $H$-Cover} problem is either polynomial-time solvable for general input graphs, or NP-complete for simple graphs on input.

\medskip
In this connection, the following relation among graphs introduced in~\cite{BFJKS21FCT} seems to play a quite important role.

\subsection{Graphs that are stronger}

\begin{definition}
A graph $A$ is said to be {\em stronger} than a graph $B$, denoted by $A\triangleright B$, if it holds true that any simple graph covers $A$ only if it also covers $B$.  
Formally, 
$$A\triangleright B\quad\Leftrightarrow\quad\forall G \mbox{ simple graph}: ((G\longrightarrow A) \Rightarrow (G\longrightarrow B)).$$
\end{definition}
We list a few of known examples and observations.

\medskip\noindent
\begin{observation}{\bf 1.}
For any two graphs $A$ and $B$, if $A\longrightarrow B$, then $A\triangleright B$.
\end{observation} 

\begin{proof}
The composition of covering projections is a covering projection. \qed
\end{proof}

\medskip\noindent
\begin{observation}{\bf 2.} 
If $A$ is a simple graph, then $A\triangleright B$ if and only if $A\longrightarrow B$ for any graph $B$. 
\end{observation}

\begin{proof}
Set $G=A$ in the definition of $A\triangleright B$. \qed
\end{proof}

\medskip\noindent
{\bf Example 1.} $F(3,0)\triangleright F(1,1)$ though $F(3,0)\not\longrightarrow F(1,1)$. A simple graph $G$ covers $F(3,0)$ if and only if it is edge 3-colorable. In such an edge-coloring by colors 1,2,3, map edges of color 1 to the semi-edge of $F(1,1)$ and the edges of colors 2 and 3 onto the loop. This is a covering projection from $G$ to $F(1,1)$.
 
\medskip\noindent
{\bf Example 2.} Let $\widetilde{P_k}$ denote the path on $k$ vertices with semi-edges incident with its terminal vertices (one semi-edge to each of them). Then $\widetilde{P_k} \triangleright C_{2k}$. A simple graph covers $\widetilde{P_k}$ if and only if it is a cycle of length divisible by $2k$, and such a cycle covers $C_{2k}$.
 
\medskip
It would certainly be too ambitious a goal to try to understand the complexity of the ``being stronger" relation, as due to Observation~2, understanding the complexity of the $\triangleright$ relation would require a full understanding of covering graphs by simple graphs, which is known to be NP-complete for many instances of the target graphs. However, there may be a hope for understanding $A\triangleright B$ for those pairs of graphs $A,B$ such that $A\not\longrightarrow B$. In the problem session of GROW 2022 workshop in Koper, September 2022, we have conjectured that the presence of semi-edges in $A$ is vital in this sense (cf.~\cite{n:GROW2022-open}).

\medskip\noindent
{\bf Conjecture 2.} If $A$ has no semi-edges, then $A\triangleright B$ if and only if $A \longrightarrow B$.

\medskip\noindent
{\bf Conjecture 3.} If $A$ has no semi-edges and no loops, then $A\triangleright B$ if and only if $A \longrightarrow B$.

\subsection{Our results}\label{subsec:ourresults}

The goal of this paper is to justify the above mentioned conjectures for several general situations. We first show that $A$ cannot be much smaller than $B$ in order to be stronger than it (in Theorem~\ref{thm:divisibility}) and we introduce two useful product-like constructions of graphs in Subsection~\ref{subsec:products}. Then we prove the conjectures for bipartite two-vertex graphs $A$ (in Theorem~\ref{thm:dipole}), for cubic graphs $A$ on at most 4 vertices, and for all (i.e., both) cubic one-vertex graphs $B$ (in Theorems~\ref{thm:3-color}~and~\ref{thm:matching}). Moreover, in the last two mentioned theorems for $B=F(3,0)$ and $B=F(1,1)$, we completely describe the graphs $A$ (even those with semi-edges) such that $A\triangleright B$. Towards this end we introduce o novel notion of \emph{graph semi-covers} in Subsection~\ref{subsec:Semicovers}. At least for cubic graphs, semi-covers fully describe the $\triangleright$ relation, as shown Theorem~\ref{thm:matching} (and also implicitly in Theorem~\ref{thm:3-color}). 

\section{General results on the $\triangleright$ relation}\label{sec:generalresults}

\subsection{Useful product-like constructions}\label{subsec:products}

We first introduce two constructions of  double covers of a given graph, one related to the direct product of this graph with $K_2$, the other one related to the cartesian product. Note here that $G^{\times}=K_2\times G$ is an often used product of $G$ and $K_2$ and is usually called the {\em universal double cover} of $G$ in the literature.

\begin{definition}
Let $G$ be a graph. 

(1) We denote by
$G^{\times}$ the graph with vertex set $V(G^{\times})=\{v_0,v_1:v\in V(G)\}$ and edge set defined as follows

(a) for every normal edge $e$ of $G$ incident with two different vertices $u, v$ in $V(G)$, there are two edges $e_0$ and $e_1$ incident with $u_0$ and $v_1$, and with $u_1$ and $v_0$, respectively, in $E(G^{\times})$;

(b) for every loop $e$ incident with a vertex $v$ of $G$, there are two parallel edges $e_0$ and $e_1$
incident with $v_0$ and $v_1$ in $E(G^{\times})$;

(c) for every semi-edge $e$ incident with a vertex $v$ of $G$, there is an edge $\widetilde e$ 
incident with $v_0$ and $v_1$ in $E(G^{\times})$.

(2) We denote by
$G^{\odot}$ the graph with vertex set $V(G^{\odot})=\{v_0,v_1:v\in V(G)\}$ and edge set defined as follows

(a) for every normal edge $e$ of $G$ incident with two different vertices $u, v$ in $V(G)$, there are two edges $e_0$ and $e_1$ incident with $u_0$ and $v_0$, and with $u_1$ and $v_1$, respectively, in $E(G^{\odot})$;

(b) for every loop $e$ incident with a vertex $v$ of $G$, there are two loops $e_0$ and $e_1$
incident with $v_0$ and with $v_1$, respectively, in $E(G^{\odot})$;

(c) for every semi-edge $e$ incident with a vertex $v$ of $G$, there is an edge $\widetilde e$ 
incident with $v_0$ and $v_1$ in $E(G^{\odot})$.
\end{definition}

\begin{figure}
\centering
{\includegraphics[width=\textwidth]{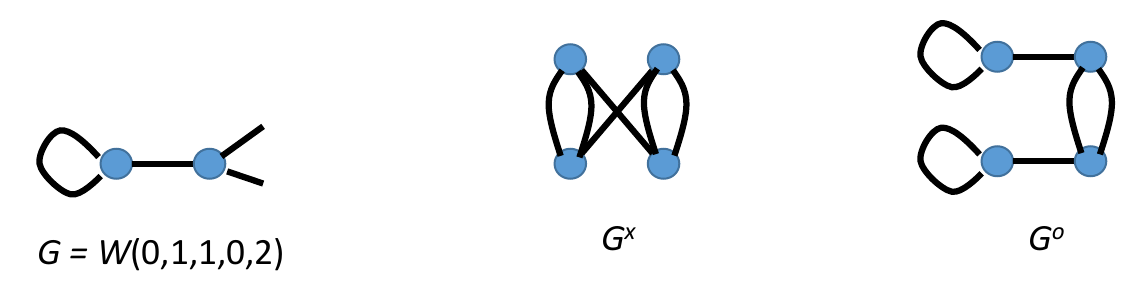}}
\caption{Examples of constructions $G^{\times}$ and $G^{\odot}$.}
\label{fig:products}
\end{figure}

See the examples in Figure~\ref{fig:products}.

\begin{proposition}\label{prop:products}
Let $G$ be a connected graph. Then both $G^{\times}$ and $G^{\odot}$ cover $G$. Moreover, $G^{\times}$ is a bipartite graph and the following hold true:

(a) $G^{\times}$ is connected iff $G$ is not bipartite, 

(b) $G^{\times}$ is the disjoint union of two copies of $G$ iff $G$ is bipartite,

(c) $G^{\times}$ is universal in the following sense: if $Y\longrightarrow G$ is a connected bipartite cover of $G$, then
$Y$ covers every connected component of $G^{\times}$.
\end{proposition}

\begin{proof} Note that by standard definition, a graph is bipartite if and only if its vertex set can be partitioned into two parts, both of them inducing edgeless subgraphs. Thus bipartite graphs contain no loops and no semi-edges.

For both $G^{\times}$ and $G^{\odot}$, the pair of mappings defined by $f_V(v_0)=f_V(v_1)=v$, $v\in V(G)$, and $f_E(e_0)=f_E(e_1)=e$, for loops and normal edges $e\in E(G)$, and $f_E(\widetilde{e})=e$ for semi-edges $e\in E(G)$, is a covering projection onto $G$.

If $A\cup B=V(G)$ is a partition of the vertex set of $G$ into two independent sets, then $G^{\times}[\{v_0:v\in A\}\cup\{v_1:v\in B\}]$ and $G^{\times}[\{v_1:v\in A\}\cup\{v_0:v\in B\}]$ are connected components of $G^{\times}$. If $G$ contains an odd cycle on vertices $u^1,u^2,\ldots,u^{2k+1}$, then $G^{\times}$ contains a cycle of length $4k+2$ on vertices $u^1_0,u^2_1,\ldots,u^{2k+1}_0,u^1_1,\ldots,u^{2k+1}_1$ and for any other vertex $u\in V(G)$, both $u_0$ and $u_1$ are connected by paths to vertices of this cycle. This proves (a) and (b).

Suppose $Y\longrightarrow G$ and let $V(Y)=A\cup B$ be a bipartition of the vertex set of $Y$. Let $f_V,f_E$ be a covering projection from $Y$ onto $G$. If $G$ is bipartite, $G^{\times}$ consists of two isomorphic copies of $G$, and $Y$ covers each of them. If $G$ is not bipartite, define mappings $g_V, g_E$ as follows
$$g_V(x)=\left\lbrace
\begin{array}{ll}
f_V(x)_0 & \mbox{ if }x\in A,\\
f_V(x)_1 & \mbox{ if }x\in B
\end{array}  
\right.
$$
and

$g_E(x)=
\widetilde{f(e)} \mbox{ if }f_E(e)\mbox{ is a semi-edge},$

$g_E(e)=f_E(e)_0$ if $e$ is incident with $x\in A$ and $y\in B$ and $f_E(e)$ is a normal edge,

$g_E(e)=f_E(e)_1$ if $e$ is incident with $x\in B$ and $y\in A$ and $f_E(e)$ is a normal edge.
 
In the case when $f_E(e)$ is a loop incident with a vertex $u\in V(G)$, $e$ belongs to an even cycle in $Y$ whose all vertices are mapped onto $u$ by $f_V$ and whose all edges are mapped onto $f_E(e)$ by $f_E$. Two-color the edges of this cycle properly and map the edges of one color onto $f_E(e)_0$ and the edges of the other color onto $f_E(e)_1$ by $g_E$. 

These mappings $g_V,g_E$ form a covering projection of $Y$ onto $G^{\times}$.
 \qed
\end{proof}

Next we show that the second construction behaves well with respect to edge-coloring. 
Since it is not so common to discuss edge-colorings of graphs with multiple edges and semi-edges, we include a formal definition here.

\begin{definition}\label{def:edge-coloring}
A {\em (proper) edge coloring} of a graph $G$ is a mapping $f:E(G)\longrightarrow C$ from the edge set of $G$ to a set $C$ of colors such that for each color $c\in C$, the subgraph $G^c=(V(G),\{e:e\in E(G), f(e)=c\}$ induced by the edges of color $c$ has maximum degree 1. The minimum size of the color set $C$ such that $G$ allows a proper coloring with these colors is called the {\em chromatic index} of $G$ and is denoted by $\chi'(G)$. We set by definition $\chi'(G)=\infty$ if $G$ has no proper edge-coloring for any number of colors (i.e., when $G$ contains at least one loop).
\end{definition}

Note that if a graph contains a loop, then it has no proper edge-coloring whatsoever.

\begin{lemma}\label{lem:indexofproduct}
For any graph $G$, $\chi'(G^{\odot})=\chi'(G)$.
\end{lemma}

\begin{proof}
If $G$ has no semi-edges, $G^{\odot}$ is the disjoint union of two copies of $G$ and the statement is clear. If $G$ has loops, so does $G^{\odot}$ and $\chi'(G^{\odot})=\chi'(G)=\infty$ since no proper edge-coloring exits, neither for $G$ nor for $G^{\odot}$.

Suppose $G$ has no loops, but does have semi-edges. If $f:E(G)\longrightarrow C$ is a proper edge-coloring of $G$ by $k=|C|$ colors, then a coloring $f^{\odot}$ obtained by keeping $f$ on the copies of normal edges of $G$ and setting $f^{\odot}(\widetilde{e})=f(e)$ for every semi-edge of $G$ is a proper $k$-edge-coloring of $G^{\odot}$. Hence $\chi'(G^{\odot})\le \chi'(G)$. On the other hand, the restriction of a proper $k$-edge-coloring $f^{\odot}$ of $G^{\odot}$ to the normal edges within $\{u_0:u\in V(G)\}$ together with setting $f(e)=f^{\odot}(\widetilde{e})$ is a proper $k$-edge-coloring of (an isomorphic copy of) $G$. Thus $\chi'(G)\le \chi'(G^{\odot})$.
\end{proof}

The concept of the $\triangleright$ relation is based on simple graphs, and the existence of simple covering graphs is guaranteed by the following theorem.

\begin{theorem}\label{thm:simplecover}
(1) Let $G$ be a graph with at most  $d$ parallel normal edges incident with the same pair of vertices and such that for every vertex, the number of semi-edges incident with this vertex plus twice the number of loops incident with it is at most $q$. Then for every even number $p\ge \mbox{max}\{d,q+1\}$, there exists a simple $p$-fold cover of $G$.

(2) If, morevoer, $G$ has no semi-edges, then a simple $p$-fold cover of $G$ exists for every $p\ge \mbox{max}\{d,q+1\}$.
\end{theorem}

\begin{proof}
Let $p$ be large enough, as specified by the assumptions. For every vertex $u\in V(G)$, take $p$ copies and denote them by $u_1,u_2,\ldots,u_p$. In a covering projection, these vertices will be mapped onto $u$.

If vertices $u\neq v\in V(G)$ are connected by $d\le p$ parallel edges, add edges of a $d$-regular bipartite graph with classes of bipartition $\{u_i:i=1,2,\ldots,p\}$ and $\{v_i:i=1,2,\ldots,p\}$. Such a bipartite graph always exists, since the complete bipartite graph $K_{p,p}$ can be properly edge-colored by $p$ colors, and we include edges of $d$ of these colors. In a covering projection, edges of one color will be mapped onto one of the parallel edges incident with $u$ and $v$ in $G$.

If $G$ has no semi-edges at a vertex, say $u$, let $\ell$ be the number of loops at $u$ in $G$. Add edges $u_iu_{i+j}$, $i=1,2,\ldots,p, j=1,2,\ldots,\ell$, addition in the subscripts being modulo $p$. This adds $2\ell$ ``vertical" edges incident with every vertex $u_i$, and in a covering projection, edges $u_iu_{i+j}$ will be mapped onto the $j$-th loop at $u\in V(G)$.

If a vertex $u$ is incident with $\ell$ loops and $s>0$ semi-edges, which means $s+2\ell\le q$, we add edges on $u_1,\ldots,u_p$ as follows. It is well known that the edges of the complete graph on $p$ vertices (note that $p$ is even in this case) can be properly colored by $p-1\ge s+2\ell$ colors. Consider such a coloring by colors $c_1,\ldots,c_{p-1}$ for the complete graph on vertices $u_1,\ldots,u_p$ and discard the edges of colors $c_i$ for $i>s+2\ell$. In a covering projection, edges of colors $c_{2i-1}$ and $c_{2i}$ will be mapped onto the $i$-th loops incident with vertex $u$, for $i=1,2,\ldots, \ell$, and edges of color $c_{2\ell+j}$ onto the $j$-th semi-edge incident with $u$, for $j=1,2,\ldots,s$. \qed   
\end{proof}

\subsection{Divisibility of the orders of graphs}\label{subsec:divisibility}

\begin{theorem}\label{thm:divisibility}
Let $A$ and $B$ be connected graphs such that  $A\triangleright B$. Then $|V(B)|$ divides $2|V(A)|$. If, moreover, $A$ has no semi-edges, then $|V(B)|$ divides $|V(A)|$.
\end{theorem}

\begin{proof}
Let $\Delta(A)$ be the maximum degree of a vertex of $A$ and let $p$ be an even number greater than $\Delta(A)$. By Theorem~\ref{thm:simplecover}, there exists a simple $p$-fold cover $G_p$ of $A$, and a simple $(p+2)$-fold cover $G_{p+2}$ as well. If $A\triangleright B$, we have $G_p\longrightarrow B$ and $G_{p+2}\longrightarrow B$. Thus $|V(B)|$ divides both $|V(G_p)|=p|V(A)|$ and $|V(G_{p+2}|=(p+2)|V(A)|$, and therefore also their difference $2|V(A)|$.

If $A$ has no semi-edges, then a simple $(p+1)$-fold cover $G_{p+1}$ of $A$ exists as well, and a similar argument straightforwardly implies that $|V(B)|$ divides $|V(G_{p+1})|-|V(G_p)|=(p+1)|V(A)|-p|V(A)|=|V(A)|$. \qed   
\end{proof}

Note that Example~2 shows a case when $A\triangleright B$ and $|V(B)|=2|V(A)|$.

\begin{theorem}\label{thm:dipole}
Let $A$ be a dipole, i.e., a graph with two vertices joined by $d$ parallel edges.
Then for every graph $B$, $A\triangleright B$ implies $A\longrightarrow B$.
\end{theorem}

\begin{proof}
By Theorem~\ref{thm:divisibility}, $B$ has at most two vertices. 

If $B$ is a one-vertex graph with
$a$ loops and $b$ semi-edges incident with the unique vertex, i.e., $B=F(b,a)$, $A\triangleright B$ implies that $d=2a+b$. Then the canonical double cover $B^{\times}$ of $B$ is isomorphic to $A$, and hence $A\longrightarrow B$.

Assume $|B|=2$. If $B$ is bipartite, then $A\cong B$ and $A\longrightarrow B$. 

Assume that $B$ is not bipartite, i. e.,
it contains loops or semi-edges. Then $B^{\times}$ is connected and has 4 vertices, and Proposition~\ref{prop:products} implies that any bipartite cover of $B$ also covers $B^{\times}$ and thus its number of vertices is divisible by 4. On the other hand, consider an odd integer $p\ge d$. By Theorem~\ref{thm:simplecover}, $A$ has a $p$-fold cover $G_p$. Such a cover is a bipartite graph and has $2p\equiv 2\mbox{ mod }4$ vertices, and thus it cannot cover $B$. Therefore $A$ is not stronger than $B$. 
\qed
\end{proof}

\subsection{A step aside: Semi-covers}\label{subsec:Semicovers}

\begin{definition}
A {\em semi-covering projection} from a graph $G$ to a connected graph $H$ is a pair of surjective mappings $f_V:V(G)\longrightarrow V(H)$ and $f_E:E(G)\longrightarrow E(H)$ such that

- $f_V$ is degree preserving,

- $f_E$ is incidence preserving (i.e., if $e\in E(G)$ is incident with vertices $u,v\in V(G)$, then $f_E(e)$ is incident with $f_V(u)$ and $f_V(v)$, which may of course be the same vertex),

- $f_E$ is a local bijection on the edge-neighborhoods of any vertex and its image.  
\end{definition}

The difference with respect to covering projections is that semi-edges may be mapped onto loops (apart from being allowed to map onto semi-edges), still loops may only be mapped onto loops, while normal edges may be mapped onto loops, semi-edges or  normal edges. The  preimage of a loop is then a disjoint union of cycles and open paths (i.e., paths with semi-edges incident to both end-vertices) spanning the preimage of its vertex, while the preimage of a normal edge is again a matching spanning the preimage of its vertex set, and the preimage of a semi-edge is a disjoint union of semi-edges and a matching spanning the preimage of its vertex.
If a graph $G$ allows a semi-covering projection onto a graph $H$, we say that $G$ {\em semi-covers} $H$, and we write $G\leadsto H$.

\begin{proposition}\label{prop:simple-semi-cover}
For any two graphs $G$ and $H$, if $G\longrightarrow H$, then $G\leadsto H$. If $G$ has no semi-edges, then $G\leadsto H$ if and only if $G\longrightarrow H$. \qed
\end{proposition}

\begin{proposition}\label{prop:composition-semi-cover}
The composition of semi-covering projections is again a semi-covering projection. \qed
\end{proposition}

\begin{proposition}\label{prop:triangle-semi-cover}
For any two graphs $A,B$, if $A\leadsto B$, then $A\triangleright B$.
\end{proposition}

\begin{proof}
If $G$ is a simple graph that covers $A$, then $G$ semi-covers $B$ by Proposition~\ref{prop:composition-semi-cover}. Since $G$ has no semi-edges, this implies that $G$ covers $B$ by Proposition~\ref{prop:simple-semi-cover}.
\end{proof}

\medskip\noindent
{\bf Example 3.} Let $A$ be a graph and let $B$ arise from $A$ by replacing two semi-edges incident with the same vertex by a loop (incident with the same vertex). Then obviously $A\leadsto B$, and hence $A\triangleright B$. This generalizes Example~1 (stating that $F(3,0)\triangleright F(1,1)$).

\section{Cubic graphs}\label{sec:Cubic}

A natural first attempt in trying to understand the $\triangleright$ relation is looking at connected regular graphs of small valency. There are only two connected 1-regular graphs, the complete graph $K_2$ and the one-vertex graph with a single semi-edge $F(1,0)$. Since $K_2\longrightarrow F(1,0)$, we also have $K_2\triangleright F(1,0)$. 

The only connected 2-regular graphs are the cycles $C_n$ and the open paths $\widetilde{P_n}$, cf. Example~2. It is easy to see that the set of connected simple covers of $C_n$ is $\{C_{kn}:k=1,2,\ldots\}$ and the set of connected simple covers of $\widetilde{P_n}$ is $\{C_{2kn}:k=1,2,\ldots\}$. Hence the following observation.

\begin{proposition}\label{prop:2-regular}
For connected 2-regular graphs, we have

- $C_n\triangleright C_m$ if and only if $m$ divides $n$;

- $\widetilde{P_n}\triangleright \widetilde{P_m}$ if and only if $m$ divides $n$;

- $C_n\triangleright \widetilde{P_m}$ if and only if $2m$ divides $n$; and

- $\widetilde{P_n}\triangleright C_m$ if and only if $m$ divides $2n$. \qed
\end{proposition}   

Note that in the first three cases of Proposition~\ref{prop:2-regular}, $A\longrightarrow B$ for the two graphs $A,B$ for which $A\triangleright B$ is stated, but in the last case, the open path $\widetilde{P_n}$ does not cover any cycle.

In the forthcoming subsections we will present results on the $\triangleright$ relation for connected cubic graphs.

\subsection{Small graphs}\label{subsec:small}

\begin{figure}
\centering
{\includegraphics[width=\textwidth]{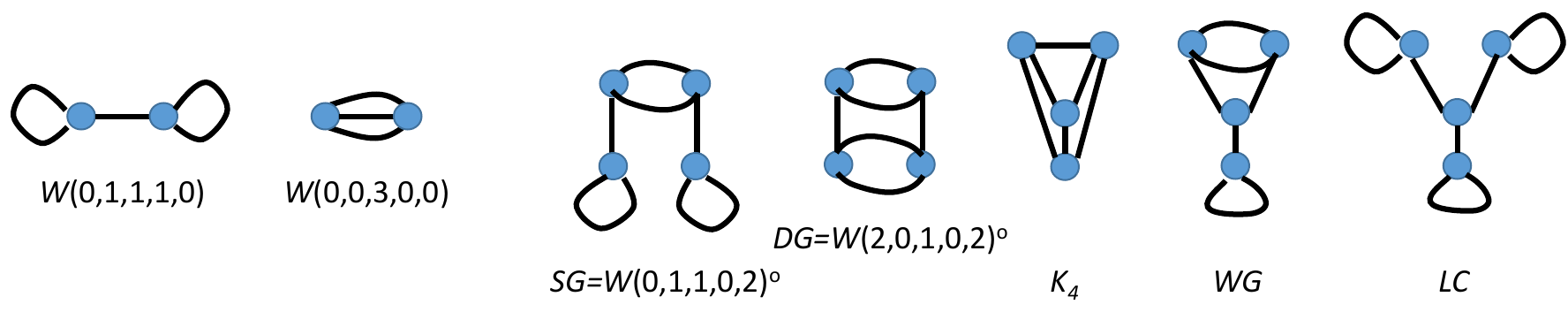}}
\caption{The connected cubic graphs with no semi-edges on 2 and 4 vertices. The notation of the 4-vertex graphs stands for Sausage Graph (SG), Drum Graph (DG), Wine Glass Graph (WG), and Loopy Claw (LC).}
\label{fig:smallcubic}
\end{figure}

There are 5 connected cubic graphs without semi-edges on 4 vertices, and 2 connected cubic graphs without semi-edges on 2 vertices (these graphs are depicted in Figure~\ref{fig:smallcubic}). We will determine all pairs $A,B$ such that $A$ is a cubic connected graph on at most 4 vertices and $A\triangleright B$. It will follow from the description that in all such cases $A\longrightarrow B$, i.e., Conjecture~2 gets justified for these graphs $A$. Theorem~\ref{thm:divisibility} implies that $|V(B)|\in\{1,2,4\}$ for every such pair $A,B$. Thus further candidates for the weaker graph $B$ are two-vertex graphs $W(0,1,1,0,2)$, $W(2,0,1,0,2)$ and $W(1,0,2,0,1)$ and one-vertex graphs $F(3,0)$ and $F(1,1)$. For a broader picture, we will include these graphs as candidates for the stronger graph $A$ as well.  

\begin{figure}
\centering
{\includegraphics[width=\textwidth]{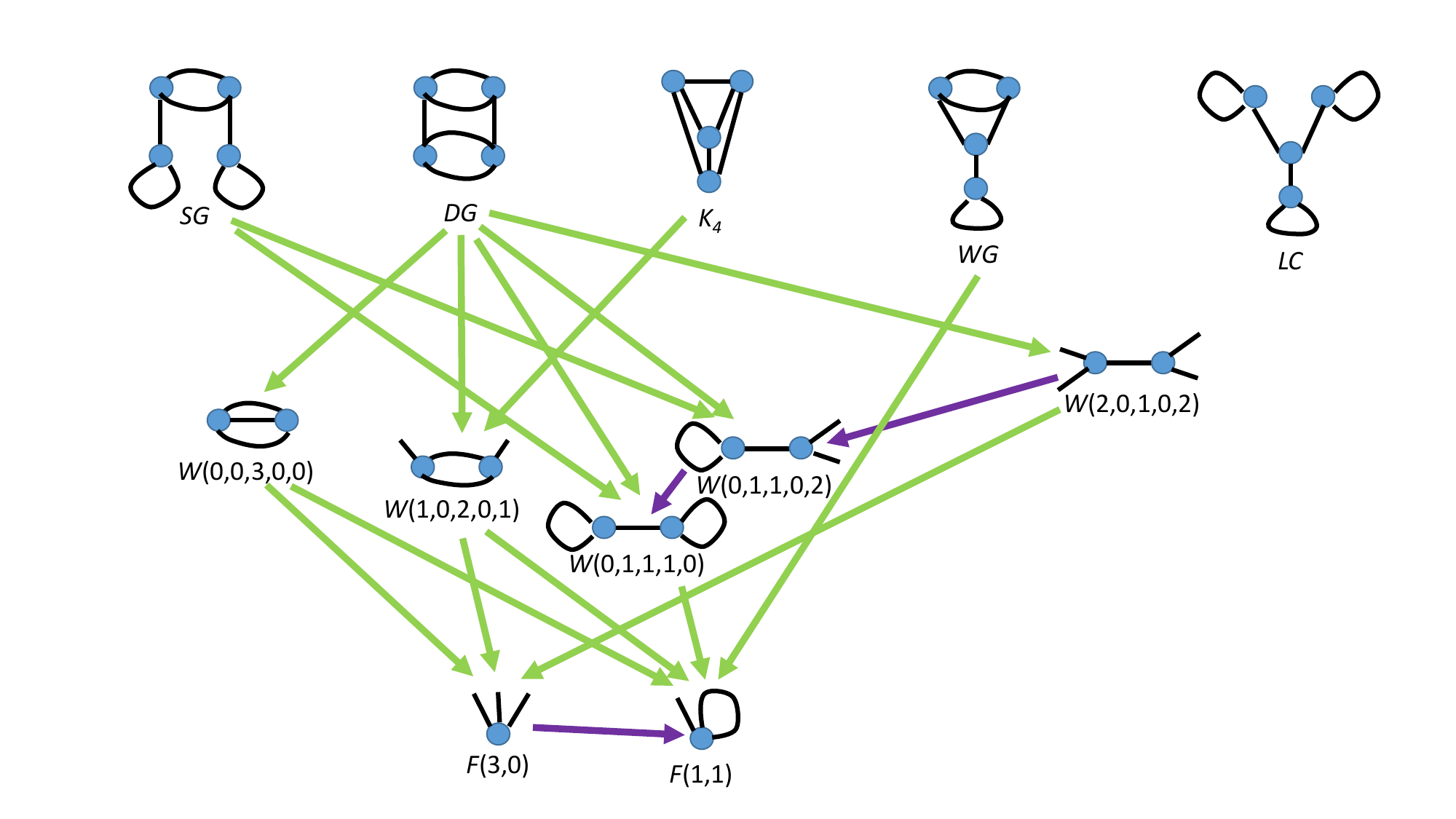}}
\caption{The partially ordered set of cubic graphs on 4, 2 and 1 vertices. The green arrows represent the existence of covering projections, the purple arrows in addition show the cases $A,B$ such that $A\triangleright B$.}
\label{fig:smallcubicposet}
\end{figure}

\begin{theorem}\label{thm:cubicsmall}
Figure~\ref{fig:smallcubicposet} shows the Hasse diagram of the partially ordered set consisting of all cubic graph with 4 vertices and no semi-edges, and all two-vertex and one-vertex cubic graphs (semi-edges allowed) ordered by existence of graph covering projections (green arrows) and the $\triangleright$ relation (purple arrows). No two graphs are in either of these relations, except those depicted in Figure~\ref{fig:smallcubicposet} and those that follow from transitivity (note that a green arrow implies a purple arrow between the same pair of graphs by   Observation~1, and we do not show these implicit purple arrows).
\end{theorem}

\begin{proof} ({\em Sketch}) It is a matter of simple checking that the green arrows express the Hasse diagram of covering projections. For the purple arrows, note that $W(2,0,1,0,2)\leadsto W(0,1,1,0,2)\leadsto W(0,1,1,1,0)$ and $F(3,0)\leadsto F(1,1)$, and hence, by Proposition~\ref{prop:triangle-semi-cover}, $W(2,0,1,0,2)\triangleright W(0,1,1,0,2)\triangleright W(0,1,1,1,0)$ and $F(3,0)\triangleright F(1,1)$.

The composition of a green and purple arrows implies the $\triangleright$ relation, and sometimes also the $\longrightarrow$ relation between the first and third graphs. As an example, consider $W(0,1,1,0,2)\triangleright W(0,1,1,1,0)\longrightarrow F(1,1)$, which implies that $W(0,1,1,0,2)\triangleright W(0,1,1,1,0)\triangleright F(1,1)$ and hence $W(0,1,1,0,2)\triangleright  F(1,1)$ by transitivity (and this purple arrow is not depicted in the figure). On the other hand, $SG\longrightarrow W(0,1,1,0,2)\triangleright W(0,1,1,1,0)$ is based on $W(0,1,1,0,2)\leadsto W(0,1,1,1,0)$ and hence $SG\leadsto W(0,1,1,1,0)$ by transitivity of the $\leadsto$ relation. But since $SG$ has no semi-edges, this implies $SG\longrightarrow W(0,1,1,1,0)$ by Proposition~\ref{prop:simple-semi-cover}. 

Showing  that no other $\triangleright$ relations hold true is more involved. For every pair $A,B$ such that $A\not\triangleright B$ is claimed, we need to provide a witness, i.e., a simple graph $G$ such that $G\longrightarrow A$ but $G\not\longrightarrow B$. This is done by a somewhat tedious case analysis, which we leave for a technical Appendix section. 
\qed
\end{proof}

\subsection{Covering $F(3,0)$}\label{subsec:snarks}

This subsection is devoted to the one-vertex graph with three semi-edges $F(3,0)$. Our aim is to show that Conjecture~2 (and therefore also Conjecture~3) hold true for this graph playing the role the weaker graph $B$. We will do so by proving the following theorem.

\begin{theorem}\label{thm:3-color}
For any graph $A$ it holds true that $A\triangleright F(3,0)$ if and only if $A\longrightarrow F(3,0)$.
\end{theorem}

Needless to say, the arguments are related to 3-edge-colorability of graphs, and to cubic graphs which are not 3-edge-colorable; such graphs are commonly called {\em snarks}. Indeed, a simple graph covers $F(3,0)$ if and only if it is 3-edge-colorable.  

\begin{proposition}
A cubic graph $G$ covers $F(3,0)$ if and only if $\chi'(G) = 3$. \qed
\end{proposition}

The following observations are simple but useful.

\begin{lemma}\label{lem:bridge}
Let $G$ be a cubic graph with no semi-edges. If $G$ contains a bridge (an edge that forms a one element edge cut), then $\chi'(G)>3$.
\end{lemma}

\begin{proof}
If the edges of $G$ were properly colored by 3 colors, edges of any two colors would occur with the same parity on every cut. But this is impossible for a cut of size one. \qed 
\end{proof}

\begin{lemma}\label{lem:cutedge}
Let $G$ be a cubic graph with no loops or semi-edges, let $e$ be an edge incident with vertices $u$ and $v$, and let $G_e$ be obtained by replacing the edge $e$ by two semi-edges $e_u,e_v$ incident with $u$ and $v$, respectively. Then $\chi'(G_e)=\chi'(G)$. 
\end{lemma}

\begin{proof}
Clearly $3\le \chi'(G_e)\le \chi'(G)\le 4$, the last inequality by Vizing theorem. Thus all we need to show is that $\chi'(G_e)=3$ implies $\chi'(G)=3$. 

If $\chi'(G_e)=3$, then $\chi'(G^{\odot})=3$ by Lemma~\ref{lem:indexofproduct}. Note that in $G^{\odot}$, the edges $\widetilde{e_u}$ and $\widetilde{e_v}$ form an edge cut of size 2. Let $f$ be a proper 3-edge-coloring of $G^{\odot}$. The parity argument implies that either every color appears on the edges $\widetilde{e_u}$ and $\widetilde{e_v}$  odd number of times, or every color appears even number of times. Since three colors cannot all occur on two edges, the former is not possible, and both edges $\widetilde{e_u}$ and $\widetilde{e_v}$ are colored by the same color. A proper 3-edge-coloring of $G$ is then obtained by using the coloring $f$ on $G-e$ and coloring the edge $e$ by the color $f(\widetilde{e_u})=f(\widetilde{e_v})$. \qed
\end{proof}

The following proposition will be useful in the proof of Theorem~\ref{thm:3-color}, but is of a more general nature.

\begin{proposition}\label{prop:bridgedsimplecover}
Let $G$ be a connected graph with no semi-edges and at least one bridge. Then $G$ has a simple cover with a bridge.
\end{proposition}

\begin{proof}
Let $e$ be a bridge of $G$ and let it be incident with vertices $u$ and $v$. Denote by $G_u$ ($G_v$) the component of $G-e$ containing $u$ (containing $v$, respectively). If $G$ has no multiple edges and no loops, then $G$ itself is a simple cover of $G$ containing a bridge. So for the rest of the proof, suppose that at least one of $G_u, G_v$ is not simple.

Let $p$ be an integer bigger than the sum of the maximum multiplicity of a normal edge in $G$ plus twice the number of loops incident to the same vertex of $G$. Then, by Theorem~\ref{thm:simplecover}, both $G_u$ and $G_v$ have simple $p$-fold and simple $(p+1)$-fold covers. We build a simple $(2p+1)$-fold cover $H$ of $G$ as follows.

\begin{figure}
\centering
{\includegraphics[width=\textwidth]{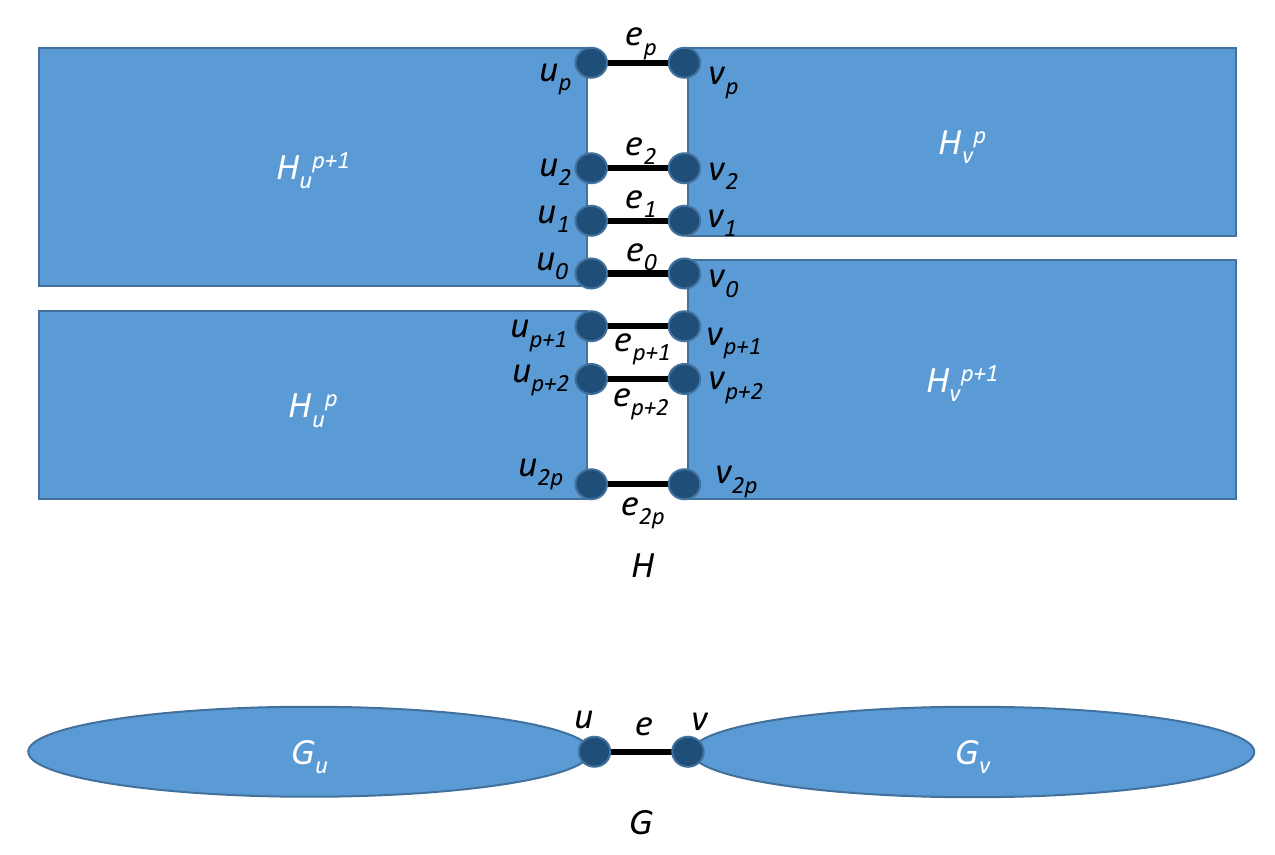}}
\caption{The construction of a simple cover with a bridge.}
\label{fig:bridgedcover}
\end{figure}

The graph $H$ will have vertices $\{x_0,x_1,\ldots,x_{2p}:x\in V(G)$. The edges are constructed so that $\{x_i:i=0,1,\ldots,2p\}$ will be the fiber of $x$, $x\in V(G)$ (in other words, the vertex mapping $f_V(x_i)=x, i=0,1,\ldots, 2p, x\in V(G)$ will be extendable to a covering projection of $H$ to $G$). Let $H_u^{p+1}$ be a simple $(p+1)$-fold cover of $G_u$ on vertices $x_i, i=0,1,\ldots,p, x\in V(G_u)$ whose existence is guaranteed by Theorem~\ref{thm:simplecover}. Similarly, let $H_u^p$ be a simple $p$-fold cover of $G_u$ on vertices $\{x_i:i=p+1,\ldots,2p, x\in V(G_u)\}$,  let $H_v^p$ be a simple $p$-fold cover of $G_v$ on vertices $\{x_i:i=1,\ldots,p, x\in V(G_v)\}$, and  let $H_v^{p+1}$ be a simple $(p+1)$-fold cover of $G_v$ on vertices $\{x_i:i=0,p+1,\ldots,2p, x\in V(G_v)\}$. Finally, we add edges $e_i$, $i=0,1,\ldots,2p$, each incident with vertices $u_i$ and $v_i$. As constructed, $H$ is a simple graph, it covers $G$ and $e_0=u_0v_0$ is a bridge of $H$. See the illustrative Figure~\ref{fig:bridgedcover}. 
\end{proof}

The following proposition is already tailored on cubic graphs.

\begin{proposition}\label{prop:k-dipolesimplecover}
Let $G$ be a connected cubic graph with no semi-edges and no loops, such that $\chi'(G)>3$. Then there exists a simple graph $H$ covering $G$ such that $\chi'(H)>3$.  
\end{proposition}

\begin{proof}
If $G$ has no multiple normal edges, then $G$ is simple and as such it itself is the desired simple cover $H$. If $G$ has a triple normal edge, then $G=W(0,0,3,0,0)$ and as such is 3-edge colorable. Suppose for the rest of the proof that $G$ has $k>0$ normal double edges and no triple one. We prove the existence of a desired graph $H$ by induction on $k$. The base case of the induction is the case $k=0$ mentioned above, where it is sufficient to take $H=G$.

\begin{figure}
\centering
{\includegraphics[width=\textwidth]{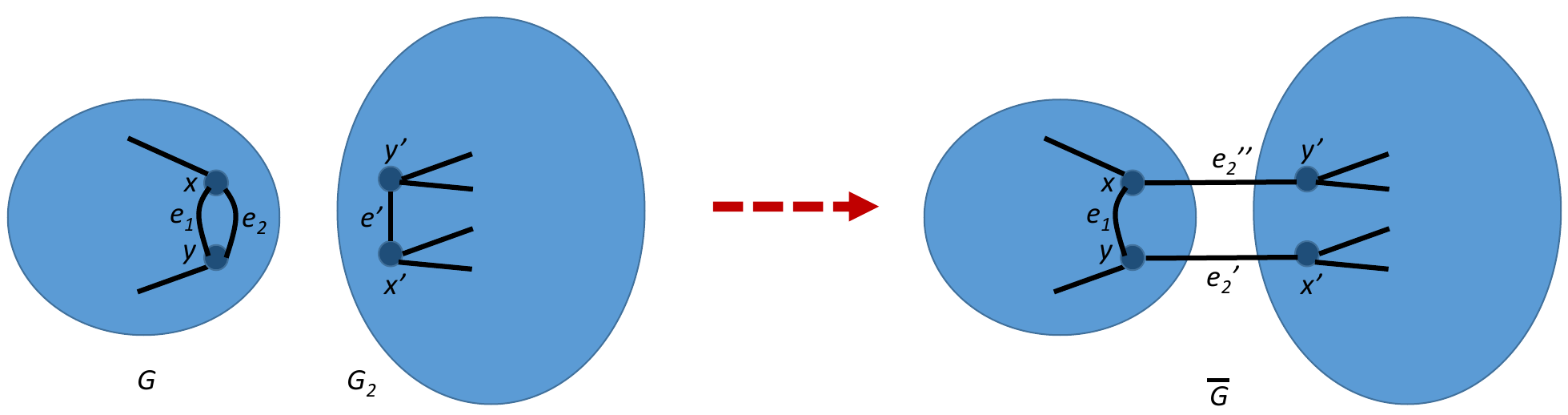}}
\caption{An illustration to the recursive construction of a simple cover.}
\label{fig:recursivesimplecover}
\end{figure}

Let $G$ contain vertices $x$ and $y$ and two normal edges $e_1,e_2$, both incident with $x$ and $y$. By Theorem~\ref{thm:simplecover}, there exists a simple connected double-cover $G_2$ of $G$. Fix a covering projection $f$ from $G_2$ to $G$. Let an edge $e'$ incident with vertices $x',y'$ of $G_2$ be such that $e'$ is mapped onto $e_2$, $x'$ is mapped onto $x$ and $y'$ is mapped onto $y$ by $f$. Construct a graph $\overline{G}$ by taking a copy of $G$ disjoint with $G_2$, deleting the edges $e_2$ (from $G$) and $e'$ (from $G_2$) and adding edges $e_2'$ connecting vertices $x'$ and $y$, and $e_2''$ connecting vertices $x$ and $y'$. Then $\overline{G}$ is a cubic graph such that

- $\overline{G}\longrightarrow G$ (take the projections induced by $f$ on $G_2$ and the identity mapping on $G$ and map both newly added edges $e_2', e_2''$ onto the edge $e_2$); 

- $\overline{G}$ is not 3-edge-colorable (if it were, any proper 3-edge-coloring would induce a proper 3-edge-coloring of $G_{e_2}$, but this is impossible since $\chi'(G_{e_2})=\chi'(G)>3$ by Lemma~\ref{lem:cutedge});

- $\overline{G}$ has $k-1$ double edges.\\
By induction hypothesis, there exists a simple graph $H$ such that $H\longrightarrow \overline{G}$ and $\chi'(H)>3$. By transitivity of the covering relation, this graph $H$ has all the desired properties. \qed  
\end{proof}

\begin{proof} ({\em of Theorem~\ref{thm:3-color}})
If $A\longrightarrow F(3,0)$ then $A\triangleright F(3,0)$ by Observation~1.

We will prove that if $A\not\longrightarrow F(3,0)$, then there exists a simple cubic graph $G\longrightarrow A$ such that $\chi'(G)>3$. Such a graph $G$ is a witness that $A\triangleright F(3,0)$ does not hold true. Hence for the rest of the proof, suppose that $\chi'(A)>3$.

\medskip
\noindent
{\em Case 1.} Assume first that $A$ has no semi-edges.

\noindent
{\em Subcase 1.1} If $A$ has a bridge, then by Proposition~\ref{prop:bridgedsimplecover}, $A$ has a simple cover $H$ which has a bridge. Thus $\chi'(H)>3$ by Lemma~\ref{lem:bridge}, and $G=H$ is the desired witness.

\noindent
{\em Subcase 1.2} If $A$ has a loop, then $A$ has a bridge, a case we have already dealt with.

\noindent
{\em Subcase 1.3} If $A$ has no loops, then by Proposition~\ref{prop:k-dipolesimplecover}, 
$A$ has a simple cover $H$ such that $\chi'(H)>3$. Thus $G=H$ is the desired witness.
  
\medskip
\noindent
{\em Case 2.} Suppose $A$ has semi-edges. By Lemma~\ref{lem:indexofproduct}, $\chi'(A^{\odot})=\chi'(A)$. Since $A^{\odot}$ has no semi-edges, the previously solved Case~1 implies that $A^{\odot}$ has a simple cover $H$ such that $\chi'(H)>3$. From the transitivity of the covering relation and from Proposition~\ref{prop:products} (stating that $A^{\odot}\longrightarrow A$), it follows that $H$ covers $A$, and thus $G=H$ is the desired witness.  
\qed
\end{proof}

\subsection{Covering $F(1,1)$}\label{subsec:matchings}

This subsection is devoted to the one-vertex ``lollipop" graph $F(1,1)$. We will not only prove that Conjecture~2 (and therefore also Conjecture~3) holds true for this graph playing the role of the weaker graph $B$, but we will also completely characterize all graphs (with any number of semi-edges) $A$ that are stronger than $B$.

\begin{theorem}\label{thm:matching}
For any graph $A$ it holds true that $A\triangleright F(1,1)$ if and only if $A\leadsto F(1,1)$.
\end{theorem}

The proof of this characterization theorem relies on several propositions and is presented in the end of this subsection. But before proving it, we show that this theorem implies the validity of Conjecture~2 for $B=F(1,1)$.

\begin{corollary}
If $A$ has no semi-edges, then $A\triangleright F(1,1)$ if and only if $A\longrightarrow F(1,1)$.
\end{corollary}

\begin{proof}
By Proposition~\ref{prop:triangle-semi-cover}, if $A$ has no semi-edges, then $A\leadsto F(1,1)$ if and only if $A\longrightarrow F(1,1)$. \qed
\end{proof}

In the rest of the subsection we will prove some auxiliary results and then present a proof of Theorem~\ref{thm:matching}. First observe that coverings of $F(1,1)$ by simple graphs are related to perfect matchings.

\begin{proposition}
A simple cubic graph covers $F(1,1)$ if and only if it contains a perfect matching. \qed
\end{proposition}

In fact all we needed in the preceding observation was that the covering graph had no semi-edges. In order to describe the general case, we introduce the following notion.

\begin{definition}
A set $M\subset E(G)$ of edges of a graph $G$ is called a {\em semi-perfect matching} if it induces a spanning 1-regular subgraph of $G$.
\end{definition} 

\begin{proposition}
A cubic graph semi-covers $F(1,1)$ if and only if it contains a semi-perfect matching. Furthermore, a cubic graph covers $F(1,1)$ if and only if it contains a semi-perfect matching which includes all semi-edges of the graph.
\end{proposition}

\begin{proof}
Let $f$ be a semi-covering projection from a cubic graph $G$ to $F(1,1)$. The preimage of the semi-edge of $F(1,1)$ is a semi-perfect matching. On the other hand, if $M\subset E(G)$ is a semi-perfect matching, its complement with respect to $E(G)$ is a spanning 2-regular subgraph of $G$. Hence mapping the edges of $M$ onto the semi-edge of $F(1,1)$ and the edges of $E(G)\setminus M$ onto the loop of $F(1,1)$ is a semi-covering projection of $G$ onto $F(1,1)$.

This semi-covering projection is a covering one if and only if $E(G)\setminus M$ is a disjoint union of cycles, which happens if and only if $M$ contains all semi-edges of $G$. \qed
\end{proof}

\begin{proposition}\label{prop:semiperfectodot}
A graph $G$ contains a semi-perfect matching if and only if $G^{\odot}$ contains a perfect matching.
\end{proposition}

\begin{proof}
If $M\subset E(G)$ is a semi-perfect matching in $G$, then $(V(G),M)^{\odot}$ is a perfect matching in $G^{\odot}$. On the other hand, if $\overline{M}$ is a perfect matching in $G^{\odot}$, then $M=\{e:e_0\in \overline{M}\}\cup \{e:\widetilde{e}\in \overline{M}\}$ is a semi-perfect matching in $G$. \qed
\end{proof}

\begin{corollary}
It can be decided in polynomial time whether an input graph with semi-edges contains a semi-perfect matching.
\end{corollary}

\begin{proof}
Given a graph $G$, we ask if $G^{\odot}$ contains a perfect matching, which can be answered in polynomial time, e.g. by Edmonds's blossom algorithm. \qed
\end{proof}

The well known Tutte theorem states that a simple graph $G$ contains a perfect matching if and only if for every subset $X\subseteq V(G)$ of vertices, the number $c_{odd}(G-X)$ of odd components of $G-X$ does not exceed the size of $X$. This theorem can be applied to graphs with multiple edges and loops, but without semi-edges. We simply discard all loops (a loop can never occur in a matching) and leave only one normal edge per a bundle of multiple edges (since at most one of them can occur in a matching), and then ask about the existence of a perfect matching in the reduced graph. For the sake of constructing a simple cover without a perfect matching, we will show that in case of cubic graphs, Tutte theorem can be stated in a stronger form. Towards this end we call a subset $X\subseteq V(G)$ {\em good} if the number of odd components of $G-X$ is strictly greater than the size of $X$, and we call $X$ {\em very good} if it is good, the subgraph $G[X]$ of $G$ induced by $X$ has no loops nor multiple edges, and no multiple edge is incident with vertices both from $X$ and from $V(G)\setminus X$ (in other words, every loop and every multiple edge is incident only with vertices from $V(G)\setminus X$).

\begin{proposition}\label{prop:strongerTutte}
Let $G$ be a connected cubic graph without semi-edges which does not contain a perfect matching. Then every inclusion-wise minimal good set of vertices is very good.
\end{proposition}

\begin{proof}
Let $X\subseteq V(G)$ be an inclusion-wise minimal good set of vertices.  By assumption, $c_{odd}(G-X)>|X|$.

If $X$ contains a vertex, say $u$, incident with a loop, such a vertex can be connected by an edge to at most 1 odd component of $G-X$. Setting $X'=X\setminus\{u\}$ we are losing at most 1 odd component, hence $c_{odd}(G-X')\ge c_{odd}(G-X)-1>|X|-1=|X'|$ and $X'$ is a good set and at the same time a strict subset of $X$, a contradiction with the assumption of minimality of $X$.

If there are two distinct vertices $u, v \in X$ joined by parallel edges in $G$, each of the vertices $u$ and $v$ can be connected by an edge to at most 1 odd component of $G-X$. Setting $X'=X\setminus\{u,v\}$ we are losing at most 2 odd components, hence $c_{odd}(G-X')\ge c_{odd}(G-X)-2>|X|-2=|X'|$ and $X'$ is a good set and at the same time a strict subset of $X$, a contradiction.

If there exist vertices $u\in X$ and $v\in V(G)\setminus X$ connected by a double edge, we set $X'=X\setminus\{u\}$ and consider several cases.

\smallskip\noindent
{\em Case 1. Vertex $u$ is adjacent to a vertex of $X$.} If the component of $G-X$ containing $v$ is odd (even), then the component of $G-X'$ containing $v$ is even (odd, respectively), no other components get changed. Hence $c_{odd}(G-X')\ge c_{odd}(G-X)-1>|X|-1=|X'|$ and $X'$ is a smaller good set, a contradiction.\\
{\em Case 2. Vertex $u$ is adjacent to a vertex in $V(G)\setminus X$ which belongs to the same component of $G-X$ as $v$.}  If the component of $G-X$ containing $v$ is odd (even), in $G-X'$ such a component is even (odd, respectively). In either case we are losing at most 1 odd component. Hence $c_{odd}(G-X')\ge c_{odd}(G-X)-1>|X|-1=|X'|$ and $X'$ is a smaller good set, a contradiction.\\
{\em Case 3. Vertex $u$ is adjacent to a vertex, say $w$, in $V(G)\setminus X$ which belongs to different component of $G-X$ than $v$.} This time we have to be more careful in the case analysis. Let us denote by $C_v$ ($C_w$) the component of $G-X$ that contains $v$ ($w$, respectively), and denote by $C'_u=C_v \cup C_w \cup \{u\}$ the component of $G-X'$ that contains $u,v$ and $w$. 

- If both $C_v$ and $C_w$ are odd, so is $C'_u$, and we lose two odd components and gain one. Hence $c_{odd}(G-X')\ge c_{odd}(G-X)-1>|X|-1=|X'|$ and $X'$ is a smaller good set, a contradiction.

- If one of $C_v, C_w$ is odd and the other one even, $C'_u$ is even and we are losing 1 odd component. Hence $c_{odd}(G-X')\ge c_{odd}(G-X)-1>|X|-1=|X'|$ and $X'$ is a smaller good set, a contradiction.

- If both $C_v$ and $C_w$ are even, $C'_u$ is odd and we are actually gaining one odd component. Hence $c_{odd}(G-X')\ge c_{odd}(G-X)+1>|X|+1>|X|-1=|X'|$ and $X'$ is a smaller good set, a contradiction. 

\smallskip
Note finally that $G$ cannot have two vertices connected by a triple edge, since then $G$ itself would be isomorphic to the dipole with triple edge, and thus $G$ would contain a perfect matching. This concludes the proof. \qed
\end{proof}

\begin{figure}
\centering
{\includegraphics[width=\textwidth]{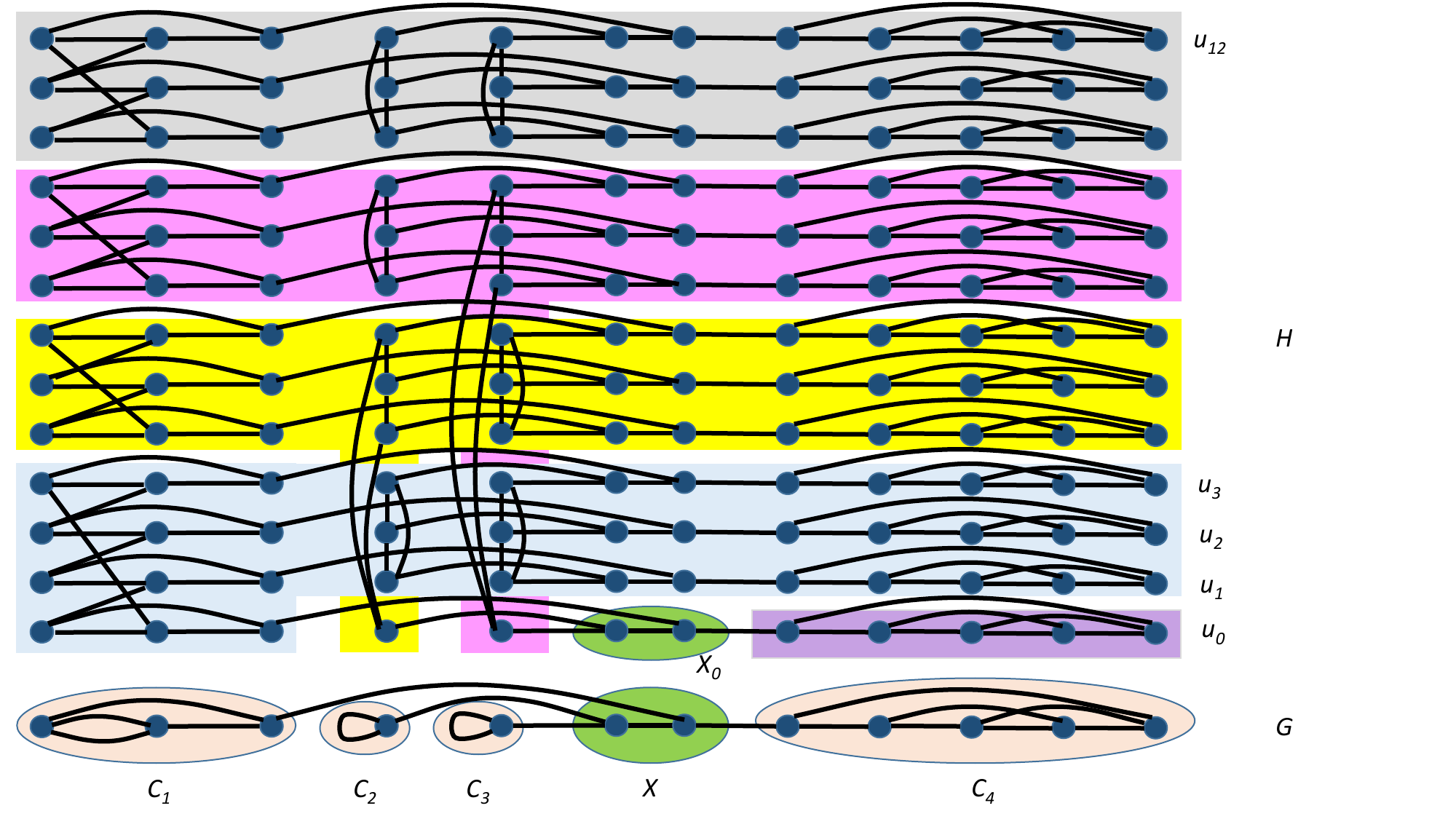}}
\caption{An illustration to the  construction of a simple cover with no perfect matching in Lemma~\ref{lem:simplecovernomatching}. The good sets ($X$ in $G$ and $X_0$ in $H$) are marked in green. The odd component $C'_1$ of $H$ is marked in light blue, the odd component $C'_2$ in yellow, the odd component $C'_3$ in pink. Since the component $C_4$ of $G$ induces a simple graph $G[C_4]$, $H[C'_4]$ is disconnected. Its even connected component is marked in grey, its odd connected component is marked in purple.}
\label{fig:Tuttecover}
\end{figure}

\begin{lemma}\label{lem:simplecovernomatching}
Let $G$ be a connected cubic graph without semi-edges which does not contain a perfect matching. Then $G$ has a connected simple cover with no perfect matching. 
\end{lemma}

\begin{proof}
If $G$ is simple, then it itself is it own simple cover without a perfect matching. Hence suppose that $G$ has a loop or a multiple edge.

By Proposition~\ref{prop:strongerTutte}, $G$ has a very good set $X$ of vertices. Let $C_1,\ldots,C_k$ be the components of $G-X$, out of which $t>|X|$ are odd. We will construct a simple $(3k+1)$-fold cover $H$ of $G$ on vertices $\{u_i:i=0,1,\ldots,3k, u\in V(G)\}$ such that $H$ does not have a perfect matching. As in the proof of Proposition~\ref{prop:bridgedsimplecover}, the projection which maps all $u_i, i=0,1,\ldots,3k$ onto $u$ for every $u\in V(G)$ can be extended to a covering projection of $H$ onto $G$.

More precisely, for every $i,j=1,2,\ldots,k$, let $H_{ij}$ be a simple cover of $G[C_j]$ (the subgraph of $G$ induced by the component $C_j$) whose existence is guarantied by Theorem~\ref{thm:simplecover}. For $i=j$, $H_{ij}$ is a 4-fold cover on vertices $\{u_0,u_{3j-2},u_{3j-1},u_{3j}:u\in C_j\}$, for $i\neq j$,   $H_{ij}$ is a 3-fold cover on vertices $\{u_{3i-2},u_{3i-1},u_{3i}:u\in C_j\}$. Further we add edges $u_iv_i$, $i=0,1,\ldots,3k$, for all normal (and simple) edges $uv$ such that at least one of the vertices $u,v$ is in $X$. Obviously, the graph $H$ constructed in this way is simple and covers $G$. 

Set $X_0=\{x_0: x\in X\}$. The components of connectivity of $H-X_0$ are subsets of $C'_j=\{u_0:u\in C_j\}\cup \{u_{3j-2},u_{3j-1},u_{3j}:u\in V(G)\}$, $j=1,2,\ldots,k$ (if $G[C_j]$ contains a loop or a multiple edge, then $H[C'_j]$ is connected and forms a component, while if $G[C_j]$ is simple, $H[C'_j]$ is the disjoint union of 3 copies of $G$ connected by edges arising from the cover of the component with a loop or a multiple edge, and one copy of $G[C_j]$). Since the total number of vertices of $G$ is even ($G$ is cubic and the hand-shaking lemma holds true even for multigraphs without semi-edges), $H$ contains exactly one odd component for every $j$ such that $C_j$ is odd. Hence $c_{odd}(H-X_0)=t>|X|=|X_0|$ and $X_0$ is a good set of vertices for $H$. Thus $H$ does not have a perfect matching. \qed  
\end{proof}

Now we are ready to present the proof of the main theorem of this subsection.

\begin{proof} ({\em of Theorem~\ref{thm:matching}}) 
If $A\leadsto F(1,1)$, then $A\triangleright F(1,1)$ by Proposition~\ref{prop:triangle-semi-cover}. For the opposite implication, we will show that for any cubic graph $A$, $A\not\leadsto F(1,1)$ implies $A\not\triangleright F(1,1)$. 

Suppose first that $A$ has no semi-edges. Then $A\not\leadsto F(1,1)$ implies that $A$ has no perfect matching. It follows from Lemma~\ref{lem:simplecovernomatching} that $A$ has a simple cover $H$ with no perfect matching. Hence $H$ does not cover $F(1,1)$ and as such is a witness that $A\triangleright F(1,1)$ does not hold true.

If $A$ has semi-edges, the assumption $A\not\leadsto F(1,1)$ implies that $A$ has no semi-perfect matching. Let us consider $A^{\odot}$. This is a  cubic graph which does not contain a perfect matching by  Proposition~\ref{prop:semiperfectodot}. Then, by Lemma~\ref{lem:simplecovernomatching}, there exists a simple cubic graph $H$ which covers $A^{\odot}$ but has no perfect matching. Transitivity of the covering relation and the fact $A^{\odot}\longrightarrow A$ imply that $H$ covers $A$, and thus it is a witness that  $A\triangleright F(1,1)$ does not hold true.
\qed
\end{proof}

\section{Acknowledgment}

Both authors gratefully acknowledge the support of Czech Science Foundation through research grant GA\v{C}R 20-15576S. We also thank Petr Hlin\v{e}n\'y for useful discussions.
 
\bibliography{bib/knizky,bib/nakryti,bib/sborniky,bib/litRN,cover}

\section*{Appendix - Full proof of Theorem~\ref{thm:cubicsmall}}

For the sake of completeness, we will show also the covering projection claims. So in the sequel, for any pair of distinct graphs $A,B$ from Figure~\ref{fig:smallcubicposet}, we will show the following

- if there is a green arrow from $A$ to $B$, then $A\longrightarrow B$,

- if there is a purple arrow from $A$ to $B$, then $A\triangleright B$ but $A\not\longrightarrow B$,

- if there is no arrow from $A$ to $B$  and neither $A\triangleright B$ nor $A\longrightarrow B$ follow from transitivity, then $A\not\triangleright B$ (and hence also $A\not\longrightarrow B$).\\
Before we process all 12 graphs from Figure~\ref{fig:smallcubicposet} in the role of the $A$ graph, we state a useful lemma. An independent set $C$ of vertices of a graph $G$ is called a {\em 1-perfect code} if every vertex from $V(G)\setminus C$ has exactly one neighbor in $C$ (such a set is also referred to as an {\em efficient dominating set} in $G$).

\begin{lemma}\label{lem:perfectcode}
If a graph $H$ contains a simple dominating vertex, then any graph $G$ that covers $H$ contains a 1-perfect code. 
\end{lemma}  

\begin{proof}
Let $u$ be a simple dominating vertex in $H$ and let $f:G\longrightarrow H$ be a covering projection. Since $u$ is a dominating vertex in $H$, its preimage $f^{-1}(u)$ is a dominating set in $G$. Since $u$ is incident with no loops or semi-edges, $f^{-1}(u)$ is an independent set. And since $u$ is incident to no multiple normal edges, every vertex $x\in V(G)\setminus f^{-1}(u)$ is dominated by exactly one vertex of $f^{-1}(u)$. \qed
\end{proof}

\medskip
\noindent
{\em Case 1. $A=F(1,1)$.} There is no arrow starting in $F(1,1)$ in Figure~\ref{fig:smallcubicposet}, thus we need to show that $F(1,1)\not\longrightarrow B$ and $F(1,1)\not\triangleright B$ for all other graphs $B$ from Figure~\ref{fig:smallcubicposet}.

For the covering projections, the only candidate for $B$ is $F(3,0)$ (no graph can cover a bigger one), but clearly $F(1,1)\not\longrightarrow F(3,0)$.

For the $\triangleright$ relation, the candidates for the graph $B$ are $F(1,1)$ and the 5 two-vertex graphs (it follows from Theorem~\ref{thm:divisibility} that every graph $B$ such that $F(1,1)\triangleright B$ has at most 2 vertices). We inspect all possible cases:

\smallskip\noindent
\underline{$1a)\quad F(1,1)\not\triangleright F(3,0)$:} By Theorem~\ref{thm:3-color}, $A\triangleright F(1,1)$ if and only if $A\longrightarrow F(3,0)$, but this is not the case for $A=F(1,1)$.\\
\underline{$1b)\quad F(1,1)\not\triangleright W(0,0,3,0,0)$:} Since $W(0,0,3,0,0)\longrightarrow F(3,0)$, $F(1,1)\triangleright W(0,0,3,0,0)$ would imply $F(1,1)\triangleright F(3,0)$, a contradiction to the previous subcase.\\
\underline{$1c)\quad F(1,1)\not\triangleright W(1,0,2,0,1)$:} Analogous to the previous subcase.\\
\underline{$1d)\quad F(1,1)\not\triangleright W(0,1,1,1,0)$:} The complete graph $K_4$ contains a perfect matching, and thus covers $F(1,1)$, but does not cover $W(0,1,1,1,0)$, since any covering projection would be a 2-fold cover and a vertex of $K_4$ mapped onto a vertex of $W(0,1,1,1,0)$ would have two neighbors mapped onto the other vertex of $W(0,1,1,1,0)$, a contradiction. So $K_4$ is a desired witness.\\
\underline{$1e)\quad F(1,1)\not\triangleright W(0,1,1,0,2)$:} Same argument as in the previous subcase, $K_4$ is a witness.\\
\underline{$1f)\quad F(1,1)\not\triangleright W(2,0,1,0,2)$:} Same argument again, $K_4$ is a witness.

\medskip\noindent
{\em Case 2. $A=F(3,0)$.} For the covering projections, the only candidate is $F(1,1)$, but $F(3,0)\not\longrightarrow F(1,1)$.

For the $\triangleright$ relation, we have $F(3,0)\leadsto F(1,1)$, and hence $F(3,0)\triangleright F(1,1)$ by Proposition~\ref{prop:triangle-semi-cover}. We will show that for every other graph $B$ from Figure~\ref{fig:smallcubicposet}, $F(3,0)\not\triangleright B$.

\smallskip\noindent
\underline{$2a)\quad F(3,0)\not\triangleright W(0,0,3,0,0), W(0,1,1,1,0), W(0,1,1,0,2), W(2,0,1,0,2)$:} The complete graph $K_4$ is 3-edge-colorable and hence covers $F(3,0)$, but the only 2-vertex graph it covers is $W(1,0,2,0,1)$. Hence $K_4$ is a witness.\\
\underline{$2b)\quad F(3,0)\not\triangleright W(1,0,2,0,1)$:} Let $K_3'$ be the complete graph with 3 vertices with 3 semi-edges added, each of them incident with a different vertex. Then $K_3'^{\odot}\longrightarrow F(3,0)$ (since $\chi'(K_3'^{\odot})=3$), but it has 6 vertices, and so it cannot cover $W(1,0,2,0,1)$ (since $W(1,0,2,0,1)$ contains semi-edges, each of its covers is even-fold). Thus $K_3'^{\odot}$ is a desired witness.  

\medskip\noindent
{\em Case 3. $A=W(0,0,3,0,0)$.} For the covering projections, $W(0,0,3,0,0)$ covers both $F(3,0)$ and $F(1,1)$, while it cannot cover any other graph from Figure~\ref{fig:smallcubicposet}.

For the $\triangleright$ relation, it follows that $W(0,0,3,0,0)\triangleright F(3,0)$ and $W(0,0,3,0,0)\triangleright F(1,1)$ and  the only other candidates are the remaining 4 two-vertex graphs (since $W(0,0,3,0,0)$ has no semi-edges, every graph $B$ such that $W(0,0,3,0,0)\triangleright B$ has at most 2 vertices by Theorem~\ref{thm:divisibility}).

\smallskip\noindent
\underline{$3a)\quad W(0,0,3,0,0)\not\triangleright W(0,1,1,1,0)$:} The complete bipartite graph $K_{3,3}$ covers the graph $W(0,0,3,0,0)$ but it does not cover $W(0,1,1,1,0)$ because the complement of any perfect matching of $K_{3,3}$ is a 6-cycle (while for a covering projection onto $W(0,1,1,1,0)$, it should be the disjoint union of 2 triangles). Thus $K_{3,3}$ is a desired witness.\\
\underline{$3b)\quad W(0,0,3,0,0)\not\triangleright W(1,0,2,0,1), W(0,1,1,0,2), W(2,0,1,0,2)$:} Again $K_{3,3}$ is a witness, since it has 6 vertices and as such cannot cover a 2-vertex graph with semi-edges (any cover of such a graph must be even-fold).

\medskip\noindent
{\em Case 4. $A=W(0,1,1,1,0)$.} This graph is a 2-fold cover of $F(1,1)$, but it is not 3-edge-colorable and hence does not cover $F(3,0)$.

For the $\triangleright$ relation, it follows that $W(0,1,1,1,0)\triangleright F(1,1)$. Since $W(0,1,1,1,0)$ has no semi-edges, every graph $B$ such that $W(0,0,3,0,0)\triangleright B$ has at most 2 vertices by Theorem~\ref{thm:divisibility}. This leaves 5 candidates.

\smallskip\noindent
\underline{$4a)\quad W(0,1,1,1,0)\not\triangleright F(3,0)$:} Since $W(0,1,1,1,0)$ is not 3-edge-colorable, it cannot be stronger than $F(3,0)$ according to Theorem~\ref{thm:3-color}.\\
\underline{$4b)\quad W(0,1,1,1,0)\not\triangleright W(0,0,3,0,0)$:} Since $W(0,0,3,0,0)$ is bipartite, it is only covered by bipartite graphs. On the other hand, $W(0,1,1,1,0)$ clearly has non-bipartite simple covers, e.g., $K_3'^{\odot}$ from Case~2 is a witness.\\
\underline{$4c)\quad W(0,1,1,1,0)\not\triangleright W(1,0,2,0,1), W(0,2,2,0,2), W(2,0,1,0,2)$:} All of these three candidate graphs contain semi-edges, and thus any simple cover is even-fold. On the other hand, $W(0,1,1,1,0)$ clearly has an odd-fold cover, e.g., $K_3'^{\odot}$ from Case~2 is a witness.  

\medskip\noindent
{\em Case 5. $A=W(1,0,2,0,1)$.} This graphs covers both $F(3,0)$ and $F(1,1)$, and hence also $W(1,0,2,0,1)\triangleright F(3,0)$ and $W(1,0,2,0,1)\triangleright F(1,1)$. Further candidates for the $\triangleright$ relation are the remaining four 2-vertex graphs and the five 4-vertex graphs, since by Theorem~\ref{thm:divisibility}, every graph $B$ such that $W(1,0,2,0,1)\triangleright B$ has at most 4 vertices.

\smallskip\noindent
\underline{$5a)\quad W(1,0,2,0,1)\not\triangleright W(0,0,3,0,0), W(0,1,1,1,0), W(0,1,1,0,2), W(2,0,1,0,2)$:} The complete graph $K_4$ covers $W(1,0,2,0,1)$, and at the same time this is the only 2-vertex graph covered by $K_4$, thanks to the symmetries of the complete graph. Thus $K_4$ is a desired witness for all four of these 2-vertex graphs.\\
\underline{$5b)\quad W(1,0,2,0,1)\not\triangleright SG, DG, WG, LC$:} Again, $K_4$ is a desired witness, since a graph covers another graph of the same order only when they are isomorphic.\\
\underline{$5c)\quad W(1,0,2,0,1)\not\triangleright K_4$:}  Let $C(n;n_1,n_2,\ldots,n_k)$ denote the cycle of length $n$ with all diagonals spanning across $n_i$ edges of the cycle, for all $i=1,2,\ldots,k$. The graph $C(8;4)$ covers $W(1,0,2,0,1)$. However, every vertex of $K_4$ is simple and dominating, but $C(8;4)$ does not contain a 1-perfect code. To see this, note that $C(8;4)$ is cubic and has 8 vertices, and thus a 1-perfect code would contain exactly 2 vertices. Because $C(8;4)$ has diameter 2, any 2 vertices are either adjacent or have a common neighbor, so no 2 vertices form a 1-perfect code. Lemma~\ref{lem:perfectcode} implies that  $C(8;4)$ does not cover $K_4$, and hence it is a desired witness.

\medskip\noindent
{\em Case 6. $A=W(0,1,1,0,2)$.} This graph is not 3-edge-colorable and thus it does not cover $F(3,0)$, and it contains a vertex incident with 2 semi-edges, and thus it does not cover $F(1,1)$. Hence it does not cover any of the graphs from Figure~\ref{fig:smallcubicposet}.

For the $\triangleright$ relation, note that $W(0,1,1,0,2)\leadsto W(0,1,1,1,0)$, and therefore also $W(0,1,1,0,2)\triangleright W(0,1,1,1,0)$. By transitivity, we get $W(0,1,1,0,2)\triangleright F(1,1)$ as well. On the other hand, $W(0,1,1,0,2)$ does not cover $F(3,0)$, and hence $W(0,1,1,0,2)\not\triangleright F(3,0)$ by Theorem~\ref{thm:3-color}.  Further candidates for the $\triangleright$ relation are the remaining three 2-vertex graphs and the five 4-vertex graphs from Figure~\ref{fig:smallcubicposet}.

\begin{figure}
\centering
{\includegraphics[width=\textwidth]{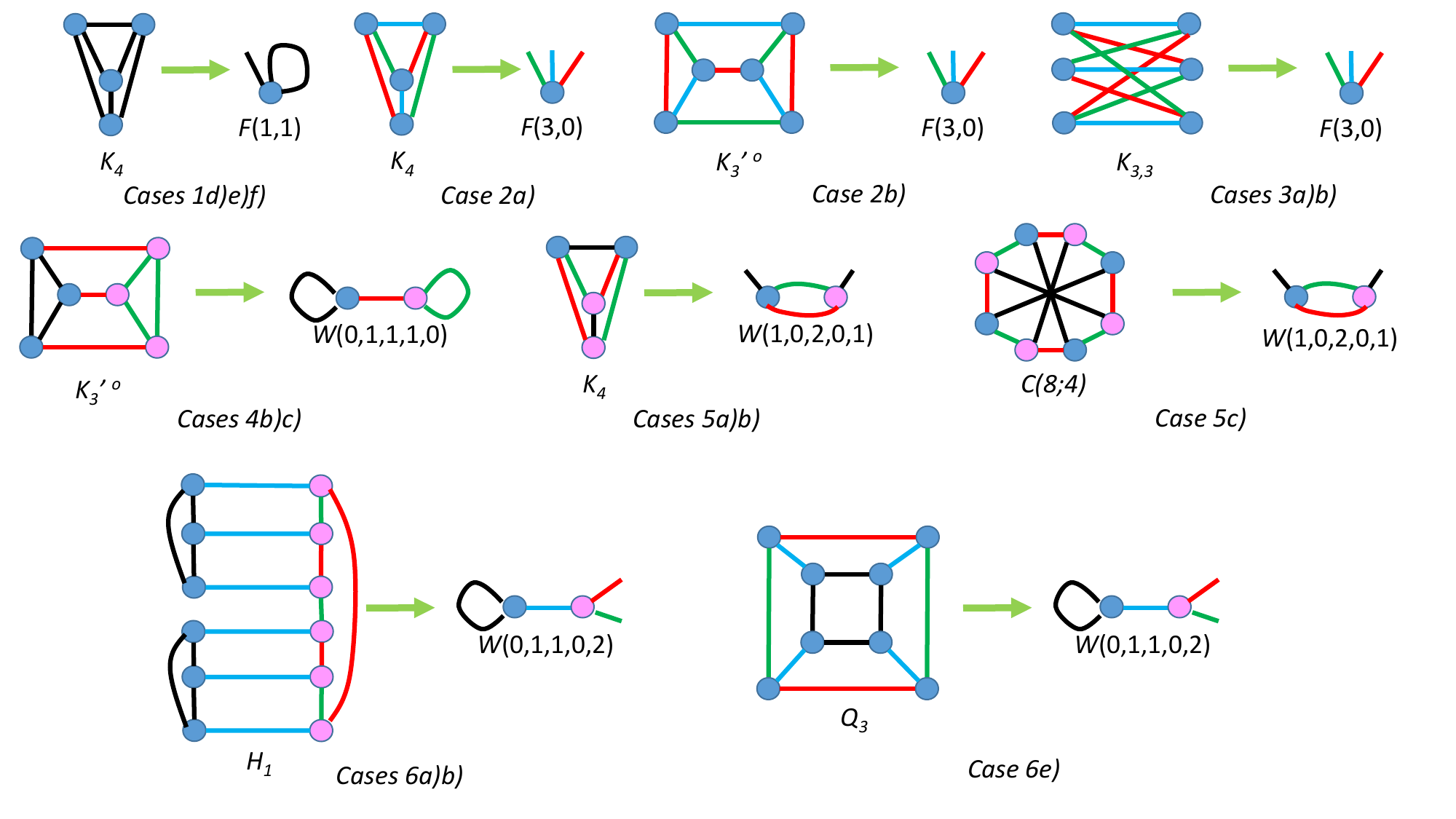}}
\caption{Witnesses for cases 1 -- 6 of the proof of Theorem~\ref{thm:cubicsmall}. Covering projections are depicted by the colors of vertices and edges.}
\label{fig:witnesses01}
\end{figure}

\smallskip\noindent
\underline{$6a)\quad W(0,1,1,0,2)\not\triangleright W(0,0,3,0,0)$:} The graph $W(0,1,1,0,2)$ has a non-bipartite simple cover (e.g., the graph $H_1$ formed by a 6-cycle on one side, the disjoint union of two 3-cycles on the other side plus a perfect matching joining the two sides), and this is a desired witness, since a non-bipartite graph cannot cover a bipartite one.\\
\underline{$6b)\quad W(0,1,1,0,2)\not\triangleright W(1,0,2,0,1)$:}
\iffalse
 In Proposition~\ref{prop:bridgedsimplecover} it was shown that a graph with semi-edges which contains a bridge has a simple cover which contains a bridge. For the proof it is sufficient that one of the components obtained after deleting the bridge has no semi-edges, and $W(0,1,1,0,2)$ has this property. Let $H$ be a simple graph with a bridge, say $uv$, such that $H\longrightarrow W(0,1,1,0,2)$. Suppose that $f:H\longrightarrow W(1,0,2,0,1)$ is a covering projection, and denote the vertices of $W(1,0,2,0,1)$ as $r$ (for {\em red}) and $b$ (for {\em blue}). Call a vertex $x$ of $H$ red if $f(x)=r$ and call it blue otherwise. The red-blue edges of $H$ (i.e., the edges with one end-point red and the other one blue) form a disjoint union of even cycles spanning the entire $H$. Thus the bridge $uv$ cannot be a red-blue edge, and so $u$ and $v$ are of the same color, say red. Consider a component of $H$ which contains $u$. Its blue vertices are paired by the matching that maps onto the semi-edge incident with vertex $b$ of $W(1,0,2,0,1)$, and all of its red vertices except for $u$ are paired by the matching that maps onto the semi-edge incident with the vertex $r$ of $W(1,0,2,0,1)$. Hence this component has an even number of blue vertices and an odd number of the red ones (including $u$). Since it has an odd number of vertices in total, it cannot be spanned by the red-blue cycles, a contradiction. Thus this $H$ is a desired witness.\\
\fi
Consider the graph $H_1$ from the previous subcase and suppose it covers $W(1,0,2,0,1)$. All 3 vertices of a triangle cannot map onto the same vertex of $W(1,0,2,0,1)$, so two of them are mapped onto one (say the blue one) and one on the other one (i.e., the pink one, if we vollow the color coding from Figure~\ref{fig:witnesses01}). But then their matched vertices on the 6-cycle must all map onto the pink vertex, and that is a contradiction.\\
\underline{$6c)\quad W(0,1,1,0,2)\not\triangleright W(2,0,1,0,2)$:} Since $W(2,0,1,0,2)$ is 3-edge-colorable, we have $W(2,0,1,0,2)\longrightarrow F(3,0)$. If it were $W(0,1,1,0,2)\triangleright W(2,0,1,0,2)$, transitivity of $\triangleright$ would imply $W(0,1,1,0,2)\triangleright F(3,0)$. However, we already know that this is not the case. Hence $W(0,1,1,0,2)\not\triangleright W(2,0,1,0,2)$ follows. \\
\underline{$6d)\quad W(0,1,1,0,2)\not\triangleright DG, K_4$:} Both $DG$ and $K_4$ cover $W(1,0,2,0,1)$. Thus if $W(0,1,1,0,2)$ were stronger than any of them, it would be also stronger than $(1,0,2,0,1)$, which we already know is not.\\
\underline{$6e)\quad W(0,1,1,0,2)\not\triangleright SG, WG, LC$:} Consider the 3-dimensional cube $Q_3$. It covers $W(0,1,1,0,2)$. If $Q_3$ covered a 4-vertex graph, it would be a 2-fold cover. The preimage of a loop in a simple graph would have at least 3 vertices, and thus $Q_3$ is not a 2-fold cover of any 4-vertex graph containing at least 1 loop. Hence $Q_3$ is a witness for all three graphs $SG, WG$ and $LC$.

\medskip\noindent
{\em Case 7. $A=W(2,0,1,0,2)$.} This graph is 3-edge-colorable and thus it covers $F(3,0)$. However, it does not cover $F(1,1)$ since it has vertices incident with 2 semi-edges each.

For the $\triangleright$ relation, we have $W(2,0,1,0,2)\leadsto W(0,1,1,0,2)$, and hence $W(2,0,1,0,2)\triangleright W(0,1,1,0,2)$ and also $W(2,0,1,0,2)\triangleright W(0,1,1,1,0)$ by transitivity of $\triangleright$. Similarly, $W(2,0,1,0,2)\leadsto F(1,1)$, and $W(2,0,1,0,2) \triangleright F(1,1)$ follows. Further candidates for the $\triangleright$ relation are $W(0,0,3,0,0), W(1,0,2,0,1)$ and the five 4-vertex graphs from Figure~\ref{fig:smallcubicposet}.

\begin{figure}
\centering
{\includegraphics[width=0.7\textwidth]{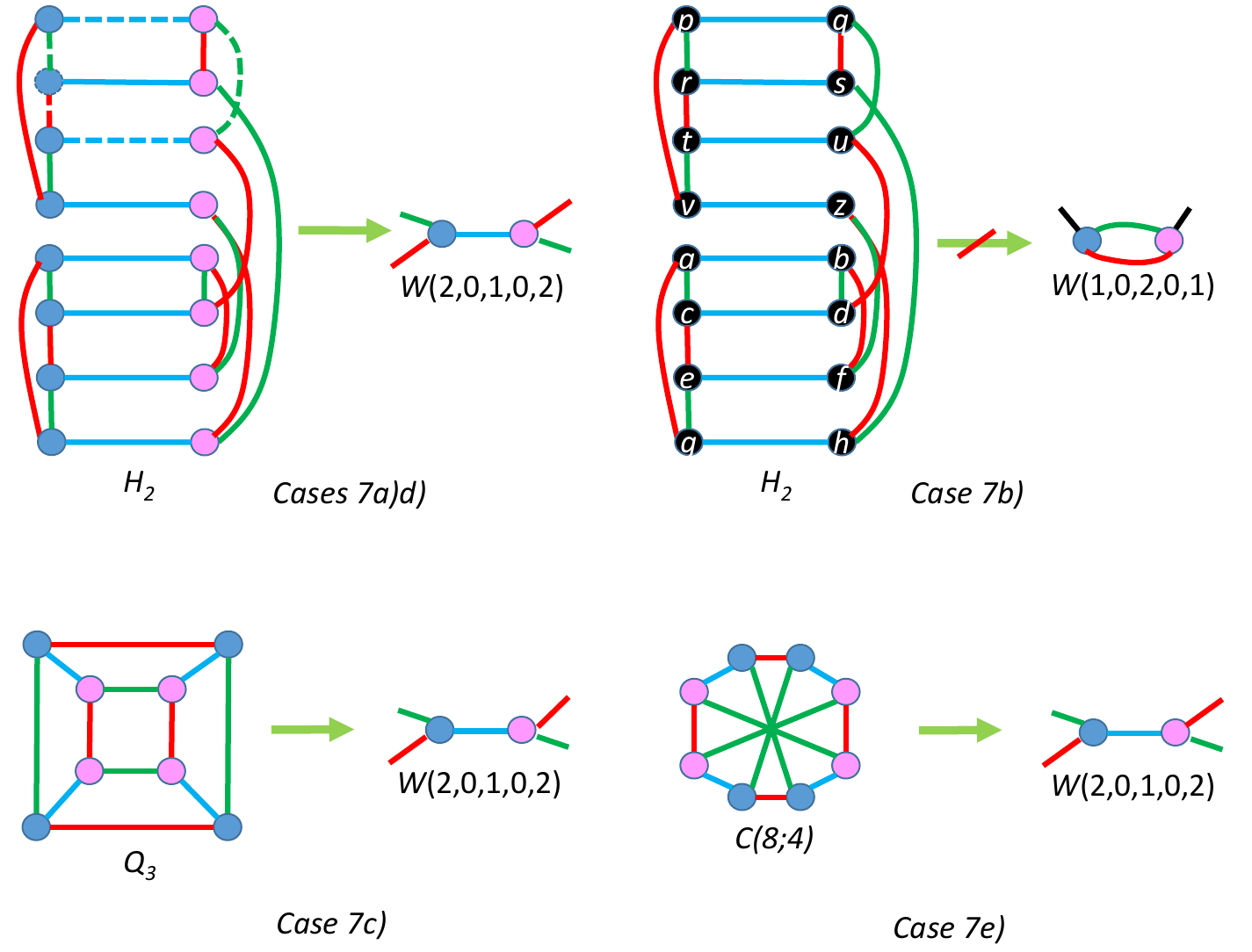}}
\caption{Witnesses for Case 7 of the proof of Theorem~\ref{thm:cubicsmall}. An odd cycle in $H_2$ is marked by dashed lines.}
\label{fig:witnesses02}
\end{figure}

\smallskip\noindent
\underline{$7a)\quad W(2,0,1,0,2)\not\triangleright W(0,0,3,0,0)$:} The graph $W(2,0,1,0,2)$ has a non-bipartite simple cover, e.g., the graph $H_2$ from Figure~\ref{fig:witnesses02} (a 5-cycle is depicted in the figure by dashed lines), and such a graph cannot cover a bipartite $W(0,0,3,0,0)$. So it is a witness.\\
\underline{$7b)\quad W(2,0,1,0,2)\not\triangleright W(1,0,2,0,1)$:} Consider the same graph $H_2$ as in the previous subcase. Suppose $H_2$ covers $W(2,0,1,0,2)$ and fix a covering projection. We say that a vertex of $H_2$ is {\em blue} if it maps onto the blue vertex of $W(2,0,1,0,2)$ in this projection, and {\em pink} if it maps onto the pink one. Note that every vertex has exactly one neighbor of its own color. It follows that every 4-cycle of $H_2$ has exactly two blue and two pink vertices. Look closer at the 4-cycle $p,q,r,s$. 

- If $p$ and $q$ are blue, and $r$ and $s$ pink, $t$ has to be blue ($r$ already has a pink neighbor $s$), $u$ has to be pink ($q$ already has a blue neighbor $p$), and $v$ has to be pink ($p$ already has a blue neighbor $q$). But then $t$ has 3 pink neighbors, a contradiction.

- If $p$ and $s$ are blue and $q$ and $r$ are pink, $t$ has to be pink (because $r$ has two blue neighbors $p$ and $s$), and hence $u$ has to be blue ($t$ already has a pink neighbor $r$). But then $q$ has 3 blue neighbors, a contradiction.

- The last possibility is that $q$ and $s$ are blue, and $p$ and $r$ are pink (or vice versa). Then $t$ has to be blue ($r$ already has a pink neighbor $p$), $u$ has to be pink ($q$ already has a blue neighbor), $v$ is blue  ($p$ already has a pink neighbor $r$), and $d$ is pink ($u$ has two blue neighbors $q$ and $t$).  Since $H_2[\{a,b,c,d,e,f,g,h\}]\equiv H_2[\{p,q,r,s,t,u,v,z\}]$, it follows that $b,e,g$ are pink and $a,c,f,h$ are blue. But then $z$ has three blue neighbors, a contradiction. So we conclude that $H_2$ is a witness for this subcase.    \\
\underline{$7c)\quad W(2,0,1,0,2)\not\triangleright SG, WG, LC$:} Similarly as in the Subcase 6e), the 3-dimensional cube $Q_3$ is a witness. It covers $W(2,0,1,0,2)$, but does not cover any of the graphs $SG, WG$ and $LC$.\\
\underline{$7d)\quad W(2,0,1,0,2)\not\triangleright DG$:} As in the Subcase 7a), a non-bipartite simple cover of $W(2,0,1,0,2)$ (e.g., $H_2$) cannot cover a bipartite graph.\\
\underline{$7e)\quad W(2,0,1,0,2)\not\triangleright K_4$:} The cycle with diagonals $C(8;4)$ covers $W(2,0,1,0,2)$ but does not cover $K_4$, as we have already seen and argued. Thus $C(8;4)$ is a desired witness. 

\medskip\noindent
{\em Case 8. $A=K_4$.} The complete graph $K_4$ is a double-cover of $W(1,0,2,0,1)$, and because of its symmetries, it does not cover any other of the 2-vertex graphs from Figure~\ref{fig:smallcubicposet}.
Then $K_4\longrightarrow F(3,0)$ and $K_4\longrightarrow F(1,1)$ follow by transitivity of the covering projection.

For the $\triangleright$ relation, it follows that $K_4\triangleright W(1,0,2,0,1)$, $K_4\triangleright F(3,0)$ and $K_4\triangleright F(1,1)$. Since $K_4$ is a simple graph, $K_4\triangleright B$ if and only if $K_4\longrightarrow B$, for any graph $B$. Thus $K_4$ is not stronger than any other of the graphs from Figure~\ref{fig:smallcubicposet}.

\medskip\noindent
{\em Case 9. $A=SG$.} The sausage graph $SG$ covers $W(0,1,1,0,2)$ because $SG=W(0,1,1,0,2)^{\odot}$. We have already argued that $W(0,1,1,0,2)\leadsto W(0,1,1,1,0)$, and hence by transitivity, $SG\leadsto W(0,1,1,1,0)$. Since $SG$ has no semi-edges, it follows that $SG\longrightarrow W(1,0,2,0,1)$, and by transitivity, $SG\longrightarrow F(1,1)$. Since the other 1- and 2-vertex graphs from Figure~\ref{fig:smallcubicposet} have no loops, $SG$ does not cover any of them.

For the $\triangleright$ relation, it follows that $SG\triangleright W(0,1,1,0,2)$, $SG\triangleright W(0,1,1,1,0)$ and $SG\triangleright F(1,1)$. For the other candidates for the $\triangleright$ relation, we argue as follows.

\begin{figure}
\centering
{\includegraphics[width=0.7\textwidth]{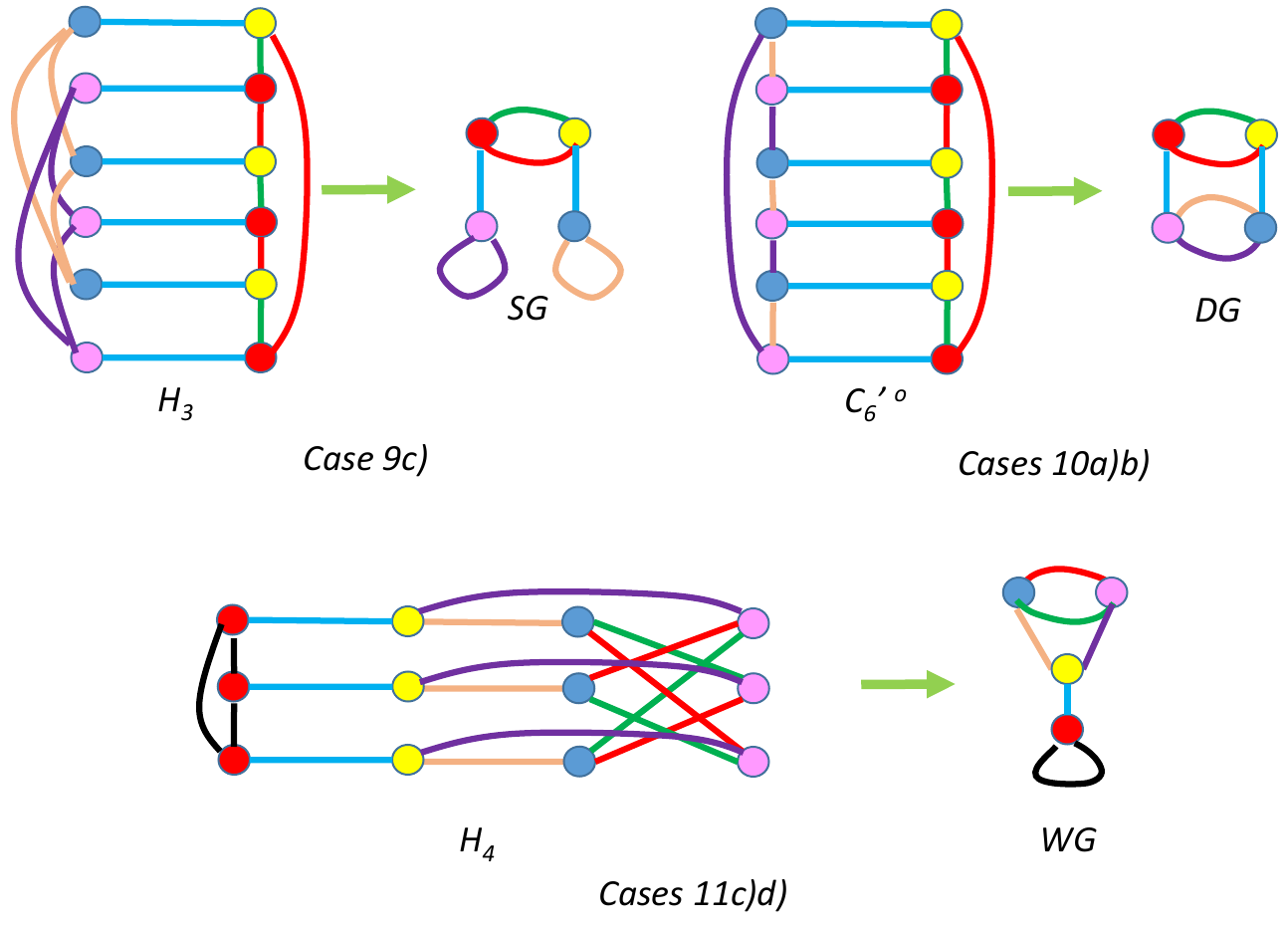}}
\caption{Witnesses for Cases 9 -- 11 of the proof of Theorem~\ref{thm:cubicsmall}. }
\label{fig:witnesses03}
\end{figure}

\smallskip\noindent
\underline{$9a)\quad SG\not\triangleright F(3,0)$:} By Theorem~\ref{thm:3-color}, $SG\mbox{$\not\triangleright$} F(3,0)$, because $SG$ is not 3-edge-colorable.\\
\underline{$9b)\quad SG\not\triangleright W(0,0,3,0,0), W(1,0,2,0,1), W(2,0,1,0,2)$:} Each of these three candidates covers $F(3,0)$, and hence if $SG$ were stronger than any of them, $SG\triangleright F(3,0)$ would follow by transitivity. However, we already know that this is not the case. Hence the claim follows.\\
\underline{$9c)\quad SG\not\triangleright WG, LC$:} The graph $H_3$ depicted in Figure~\ref{fig:witnesses03} covers $SG$, but does not contain a 1-perfect code, and hence cannot cover either of the graphs $K_4, WG, LC$. A 1-perfect code in $H_3$ would have 3 vertices, at most one in each of the triangles, and at most 2 in the 6-cycle. Hence it contains at least one vertex in the 6-cycle, and the triangle containing the neighbor of this vertex cannot contain any code-vertex. It follows that each vertex of this triangle is dominated by a different code-vertex from the 6-cycle, but that is impossible.

\medskip\noindent
{\em Case 10. $A=DG$.} The drum graph covers all 2-vertex (and also all 1-vertex) graphs from the figure. Hence it is also stronger than all of them. We will argue about the 4-vertex candidates.

\smallskip\noindent

\underline{$10a)\quad DG\not\triangleright K_4, WG, LC$:} 
 Let $C_6'$ denote the cycle of length 6 with a semi-edge attached to every vertex. Then the graph $C_6'^{\odot}$ (a graph often called the {\em prism over $C_6$}) does not have a 1-perfect-code (a 1-perfect code would have 3 vertices, 2 of them in one of the 6-cycles, in which case the undominated vertices of the other 6-cycle would induce a path, not a star, and thus could not be dominated by a single vertex). Therefore, by Lemma~\ref{lem:perfectcode}, $C'^{\odot}$ covers neither $K_4$, nor $WG$, nor $LC$. Since $DG\longrightarrow C'^{\odot}$, it is a desired witness.\\
\underline{$10b)\quad DG\not\triangleright SG$:} The prism graph $C_6'^{\odot}$ is again a desired witness. If it covered $SG$, it would be a 3-fold cover and every the preimage of every vertex of $SG$ would have size 3. Thus the preimage of a vertex incident with a loop would necessarily induce a 3-cycle in the covering graph, but $C'^{\odot}$ is triangle-free.

\medskip\noindent
{\em Case 11. $A=WG$.}
The wine glass graph covers $F(1,1)$ and no other graph from Figure~\ref{fig:smallcubicposet}. 

Thus we have $WG\triangleright F(1,1)$. We will argue that $F(1,1)$ is the only graph for which $WG$ is stronger.

\smallskip\noindent
\underline{$11a)\quad WG\not\triangleright F(3,0)$:} By Theorem~\ref{thm:3-color} only (and exactly) 3-edge-colorable graphs are stronger than $F(3,0)$, and $WG$ is not 3-edge-colorable.\\
\underline{$11b)\quad WG\not\triangleright W(0,0,3,0,0), W(1,0,2,0,1), W(2,0,1,0,2)$:} Each of these graphs covers $F(3,0)$, and thus $WG\triangleright B$ for any $B\in \{W(0,0,3,0,0), W(1,0,2,0,1), W(2,0,1,0,2)\}$ would imply $WG\triangleright F(3,0)$ what we already know is not the case.\\
\underline{$11c)\quad WG\not\triangleright W(0,1,1,1,0)$:} Consider the graph $H_4$ from Figure~\ref{fig:witnesses03}. It covers the wine glass graph. Suppose it covers $W(0,1,1,1,0)$ and fix a covering projection. Call the vertices of $W(0,1,1,1,0)$ blue and pink. All three vertices of the triangle must map onto the same vertex of $W(0,1,1,1,0)$, say on the blue one. Then their neighbors in $H_4$ must map onto the pink vertex of $W(0,1,1,1,0)$. Each of these neighbors has already one neighbor mapped on the blue vertex, and so the other two neighbors map onto the pink vertex as well. But that means that 3 vertices of $H_4$ are mapped onto the blue vertex and 9 of them onto the pink one, and the mapping cannot be a covering projection.\\
\underline{$11d)\quad WG\not\triangleright W(0,1,1,0,2)$:} We already know that $WG\mbox{$\not\triangleright$} W(0,1,1,1,0)$. Hence $WG\triangleright W(0,1,1,0,2)$ would lead to a contradiction, namely $WG\triangleright W(0,1,1,1,0)$.\\
\underline{$11e)\quad WG\not\triangleright SG, DG$:} Each of these graphs covers $W(0,1,1,0,2)$, and hence both $WG\triangleright SG$ and $WG\triangleright DG$ would lead to a contradiction, namely $WG\triangleright W(0,1,1,0,2)$ what, as we already know, does not hold true.\\
\underline{$11f)\quad WG\not\triangleright K_4$:} The complete graph $K_4$ covers $F(3,0)$, and hence $WG\triangleright K_4$ would imply $WG\triangleright F(3,0)$, a contradiction.\\
\underline{$11g)\quad WG\not\triangleright LC$:}Consider again the graph $H_4$ from Subcase 11c). Any covering projection onto $LC$ would be a 3-fold cover, and thus the preimages of the 3 vertices of $LC$ incident with loops would induce triangles in the covering graph. However, $H_4$ has only one triangle. Thus $H_4$ is a witness.

\medskip\noindent
{\em Case 12. $A=LC$.} The loopy claw is neither 3-edge-colorable, nor does it contain a perfect matching. Hence it neither covers nor is stronger than any other of the graphs from Figure~\ref{fig:smallcubicposet}.
   
\qed

\end{document}